\documentclass[11pt,twoside,english]{article}
\usepackage{
babel,graphicx,psfrag,amsfonts}
\newcommand{\fra}[2]{\displaystyle\frac{\mathstrut#1}{\mathstrut#2}}

\newcounter{cuent}
\newcounter{theorem}[section]
\newtheorem{defi}[theorem]{\sc Definition}
\newtheorem{lema}[theorem]{\sc Lemma}
\newtheorem{prop}[theorem]{\sc Proposition}
\newtheorem{cor}[theorem]{\sc Corollary}
\newtheorem{teoremadeotro}[theorem]{\sc Theorem}
\newtheorem{teo}{\sc Theorem}

\newtheorem{obs}[theorem]{\sc Remark}
\newtheorem{parag}[theorem]{\sc}

\newcommand{\A}[6]{#1.  {\em #2. \em} #3 \mbox{\bf #4} (#5),
pp #6.}
\newcommand{\B}[4]{#1. {\em #2. \em} #3 (#4).}

\begin{document}

\title{Simultaneous Continuation of Infinitely Many  Sinks Near a Quadratic Homoclinic Tangency.}
\author{Eleonora Catsigeras\thanks{  E-mail: eleonora@fing.edu.uy   }, \  Marcelo Cerminara\thanks {E-mail: cerminar@fing.edu.uy} \ and Heber Enrich\thanks{
 E-mail: enrich@fing.edu.uy
   Instituto de Matem{\'a}tica y Estad{\'\i}stica Rafael Laguardia (IMERL),
 Fac. Ingenieria. Universidad de la Rep{\'u}blica.  Uruguay. Address: Herrera
  y Reissig 565. Montevideo. Uruguay. }}

\date{October 27th, 2008.}
\maketitle

\begin{abstract}

We prove that the 
$C^3$ diffeomorphisms on surfaces, exhibiting infinitely many sinks
near the generic unfolding of a quadratic homoclinic tangency of a
dissipative saddle, can be perturbed along an infinite dimensional
manifold of $C^3$ diffeomorphisms such that infinitely many sinks
persist simultaneously. On the other hand, if they are perturbed
along one-parameter families that unfold generically the quadratic
tangencies, then at most a finite number of those sinks have
continuation.

\end{abstract}

\section{Introduction and statement of the main results.} \label{section1}

Let $M $ be a two-dimensional $C^{\infty}$ compact and connected
riemannian manifold, and let $\mbox{Diff}^3(M)$ be the infinite
dimensional manifold of all $C^3$-diffeomorphisms $f \colon M
\mapsto M$.

Let  ${f}_0 \in \mbox{Diff}^3(M)$ having a saddle fixed point $P_0$.
We denote ${\lambda}_0 < 1 < {\sigma}_0$ the eigenvalues of $D{f}_0
(P_0)$.

We consider  diffeomorphisms that are dissipative in a saddle
point, i.e. ${\lambda}_0 {\sigma}_0 < 1$.
We also  assume  that the diffeomorphism ${f}_0 $  exhibits at $q_0
$ a quadratic homoclinic tangency (see \cite{palistakens}) of the
saddle point $P_0$, recalling the following definition:

\begin{defi} \em \label{quadratic}

 We say that the homoclinic tangency at $q_0$ of the periodic saddle
point $P_0$ is quadratic if there exists a $C^2$ local chart in a neighborhood  of $q_0$ such that the
 stable arc of $P_0 $  which contains the tangency point $q_0 $ has equation $y=0$,
and the  unstable arc has equation  $y={\overline{\beta}}
\, x^2$ with $ {\overline{\beta}}\neq 0$. \em
\end{defi}

Take a one-parameter family $\{\widetilde{f}_t\}_{t \in I} \subset
\mbox{Diff}^3(M)$ through the given map $\widetilde{f}_0=f_0$, such that the
quadratic homoclinic tangency unfolds generically  into two
transversal homoclinic intersections for $t >0$.

The Newhouse-Robinson Theorem (\cite{newhouse}, \cite{robinson})
asserts that, as near as wanted from ${\widetilde{f}}_0$ in the one-parameter
family $\{\widetilde{f}_t\}_{t\in I}$, there exists an interval $I_0$ and a dense set
$J_0\subset I_0$ of values of the parameter such that for all $t \in J_0$,
$\widetilde{f}_t$ exhibits infinitely many simultaneous sinks.

We will prove that for values $t$ in a dense set  $J\subset J_0$,
the map $\widetilde{f}_t$ is bifurcating: in fact, our Theorem
{\ref{teoprevio}} asserts that at most a finite number of certain sequence of infinitely many sinks
of $\widetilde{f}_t$ can simultaneously persist when we perturb
$\widetilde{f}_t$ along certain one-parameter families in
${\mbox{Diff}^3(M)}$. Nevertheless, in Theorem \ref{teo} we prove
that  the bifurcation of infinite many
simultaneous sinks has infinite dimension in Diff$^3(M)$. 


Now let us define the kind of perturbations of each diffeomorphism  and the kind of persistence of each sink which we will consider all along this paper:

%
%
%

\begin{defi} \em

Let us suppose that $g_0 \in \mbox{Diff}^3(M)$  exhibits a sink
$s_0$. Consider $g_1 \in \mbox{Diff}^3(M)$, isotopic to
$g_0$. We say that $g_1$ exhibits the continuation $s_1 = s(g_1)$ of
the sink $s_0$, if there is a differentiable isotopy $\{g_t\} _{t
\in \mathbb{R}} \subset {\mbox{Diff}^3(M)}$ such that for all $t \in
[0,1]$ there exists a sink $s_t = s(g_t)$ of $g_t$ and the
transformation $t \in [0,1] \mapsto s_t \in M$ is of $C^1$ class.

\end{defi}

We are now ready to state the main result of this paper:


{\teo \label{teo}

Let $M$ be  a  $C^\infty$ two dimensional compact connected
riemannian manifold. Let  \em  ${f}_0 \in \mbox{Diff}^3(M)$ \em
exhibiting a quadratic homoclinic tangency of the saddle point
$P_0 $. Assume  that the saddle is dissipative, i.e. its eigenvalues
${\lambda}_0 < 1 < {\sigma}_0$
  verify
${\lambda}_0 {\sigma}_0 <1$.

Then, given an arbitrarily small neighborhood \em ${\cal N}$ of
${f}_0$ in $\mbox{Diff}^3(M)$ \em
there exists a $C^1$
arc-connected infinite-dimensional local submanifold $\mathcal{M}
\subset {\cal N}$  such that:

(a) Every $g\in
\mathcal{M}$ exhibits infinitely many simultaneous sinks $s_i(g)_{i \in \mathbb{N}}$.

(b) Each sink $s_i(g)$ is the   continuation of the respective sink
$s_i(g_0)$, for any pair of diffeomorphisms
  $g_0, g \in \mathcal {M}$.}

\vspace{.2cm}

 Note that the given diffeomorphism ${f}_0$
does not necessarily belong to ${\cal M}$. We prove theorem 1
through sections 2 to 8.



\begin{obs} {Notation} \label{defincionV} \em

Let  $f_{} \in \mbox{Diff}^3(M)$ have a horseshoe $\Lambda \subset
M$, as defined in \cite{palistakens} Chapter II, Section 3. As
$\Lambda$ is an hyperbolic set, there exist   constants $C >0$,
$\widetilde\lambda < 1 $ and $\widetilde\sigma> 1$ and a splitting
$T_pM= E^u \oplus E^s $ for all $p \in \Lambda$, such that $||
Df_{}^n(v)||\geq C \widetilde\sigma^n ||v|| \;\;\forall v \in E^u$
and $|| Df_{}^n(v)||\leq C^{-1} \widetilde\lambda^n ||v|| \;\;
\forall v \in E^s$. Besides  the horseshoe $\Lambda$ is the maximal
invariant set in an open neighborhood ${U}$ of itself. In Section
\ref{seccionCoordenadas} we will add some other restriction to $U$.

{\begin{figure}
[h]\psfrag{omega}{$\Omega$}\psfrag{v}{$V$}\psfrag{q}{$q$}
\psfrag{u}{${U}$} \psfrag{s}{$\sigma$}\psfrag{l}{$\lambda$}
\psfrag{fn}{$f^n$}\psfrag{fN}{$f^N$}
\begin{center}\includegraphics[scale=.5]{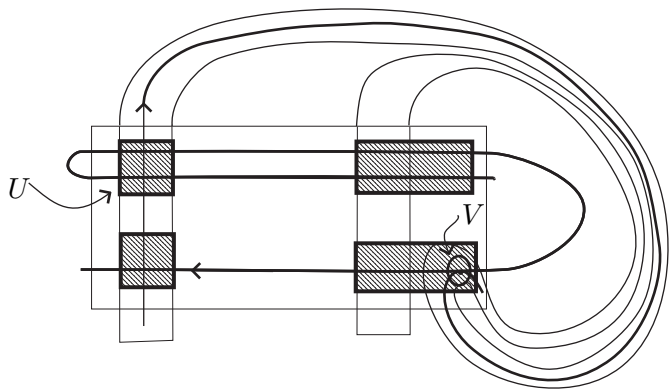}
\caption{\label{gr100} Horseshoe}\end{center}
\end{figure}}


We assume that $f_{}$ exhibits at the point $q_1 \in U$ a quadratic
homoclinic tangency of the invariant manifolds of a periodic saddle
point $P_1 \in \Lambda$. We choose $q_1$ such that $f_{}^{-1}(q_1)
\not \in U$. It is not restrictive to consider $q_1 \in W^s(P_1)$
such that for all $n\geq 1$, $f_{}^n(q_1) $ belongs to $ U$ and then
call $N_1
>1$ to an integer number such that $f_{}^{-N_1}(q_1) \in U$
belongs to the local stable manifold of $P_1$. We will take a small
neighborhood $V$ of $q_1$ such that $V \subset U$, $f_{}^{-N_1}(V)
\subset U$, and $f_{}^{-i}(V)$ and $\Lambda$ are pairwise disjoint
for $i=0,\ldots N_1$. We also assume that $f_{}^{-1}(V) \cap U
=\emptyset$. We will work with small perturbations of $f_{}$, such
that the former properties of $V$ persist.

As shown in \cite{palistakens}, Chapter II Section 3 and Appendix 1,
it is possible to construct invariant stable and unstable local
foliations $\mathcal{F}^s$ and $\mathcal{F}^u$ in a neighborhood of
$\Lambda$. We will denote $W^s_{\mbox{\footnotesize loc}}$ and
$W^u_{\mbox{\footnotesize loc}}$ the respective leaves of the
foliations. These foliations are $C^{1+\varepsilon}$, meaning in
particular that the tangent directions to the leaves are
$C^{1+\varepsilon}$.
\end{obs}

\begin{defi} \label{df}
\em The line of tangencies $L(f_{})$ is the set of points in a small
neighborhood $V$ of $q_1$ where the leaves of $\mathcal{F}^s$ and
$\mathcal{F}^u$ are tangent. \em
\end{defi}

\begin{obs}
\label{tangenciascuadraticastodas} \em Since the tangent directions
are $C^{1+\varepsilon}$, the tangencies on $L(f_{})$ are also
quadratic, and $L(f_{})$ is a differentiable curve (see
\cite{palistakens} Chapter V, Section 1). It persists and depends
continuously on $f_{}$. \em
\end{obs}

%
%
%
%
%
%
\begin{defi}
\label{Cantorsets} \em The stable (unstable) Cantor set $K^s$ (resp.
$K^u$) is  the intersection with the line of tangencies $L(f_{})$ of
the local leaves through the points $P \in \Lambda$ of the stable
foliation $\mathcal{F}^s$  (resp. the $f_{}^{N_1}$ iterates of the
local leaves of $\mathcal{F}^u$ passing through the points $P \in
\Lambda$).
\end{defi}

\begin{defi} \em \label{defiunfolding}
We say that a one-parameter family
$\{\widetilde{f}_t\}_{-\varepsilon\leq t \leq + \varepsilon} \subset
{\mbox{Diff}^3(M)}$ passing through
a 
diffeomorphism $f_{0}$, unfolds generically the  quadratic
tangencies of the horseshoe $\Lambda (f_{0})$, if there exists a
velocity $v
>0$ such that
$$\left|\frac{d \mu _{p,q}(\widetilde{f}_t)} {dt} \right| \geq v >0 \; \; \forall p,q
 \in \Lambda, \; \forall t \in (-\varepsilon,   + \varepsilon)$$
where $\mu_{p, q}(\widetilde{f}_t)$ is the distance along the line
of tangencies $L(\widetilde{f}_t)$ between $p^s  =   $
$W^s_{\mbox{\footnotesize loc}}(p) \bigcap L(\widetilde{f}_t) \in
K^s$  and $q^u \in f_t^{N_1} (W^u_{\mbox{\footnotesize loc}} (q))
\bigcap L(\widetilde{f}_t) \in K^u$ of any two points $p, q \in
\Lambda (\widetilde{f}_t)$.
\end{defi}

\begin{teo}
\label{teoprevio} In the hypothesis of Theorem \em \ref{teo}, \em
given a one-parameter family $\{\widetilde{f}_t\}_{-\varepsilon \leq
t \leq \varepsilon} \subset \mbox{Diff}^3(M)$ which generically
unfolds the quadratic homoclinic tangency at $q_0$ exhibited by
$f_0$, there exist an open real interval $I\subset
(-\varepsilon,\varepsilon)$ and a dense set $J \subset I$ of the
parameter values such that
if $f_\infty \in \{\widetilde{f}_t\}_{t \in J}$,  then:

 \vspace{.2cm}
\begin{description}

\item{(A)} \label{q} $f_\infty$ exhibits infinitely many sinks $s_i(f_\infty)_{i \geq 1} \in V$
with  periods $p_i (\rightarrow +\infty)$, and there exists $0<\rho
<1$ such that the eigenvalues of $df_\infty^{p_i}(s_i)$ have modulus
smaller than $\rho$ for all $i \geq 1$.

\item{(B)}
There exists a local $C^1$ infinite-dimensional, arc-connected
manifold \em ${\cal M} \subset {\mbox{Diff}^3 (M)}$\em,
 such that:

\begin{enumerate}

\item \label{i} $f_\infty \in {\cal M}$

\item \label{ii} If $g \in {\cal M}$ then $g$ exhibits the  continuation
$s_{i}(g) \in V$ of the infinitely many  sinks $s_{i}(f_\infty)$.
\end{enumerate}

%
%

\item{(C)} Any one-parameter family $\{g_\mu\}_{-\varepsilon\leq \mu \leq + \varepsilon}$
of $C^3$ diffeomorphisms passing through $g_0=f_\infty$ and
unfolding generically the    quadratic tangencies on $L(g_0)$,
exhibits for $\mu \neq 0$ at most a finite number of simultaneous
continuations $s_i(g_\mu)$ of the sinks $s_i(g_0)$ constructed in
part (A).
\end{description}

\end{teo}
Remark to thesis (A): The sinks $s_i (\widetilde{f}_t)$ for $t \in J$ are not
necessarily the continuation of the sinks $s_i(f_\infty)$, at least not
for infinitely many values of $i \geq 1$.

We prove Theorem \ref{teoprevio} in Section \ref{seccionpruebateo}.


\vspace{.2cm}

An interesting open problem is the prevalence of infinite sinks. A
conjecture of Palis (\cite{palis}) asserts that there exists a dense
set of $C^r$ diffeomorphisms with a finite number of attractors with
a total Lebesgue measure attracting basins. In dimension two the
main obstruction to this conjecture is that the phenomenon of the
coexistence of infinite simultaneous sinks occurs for a whole open
set in $\mbox{Diff}^r(M)$. We observe that our
Theorem 2 does not solve the problem, since infinite sinks could
appear from other homoclinic tangencies.

To prove Theorem \ref{teo}, inspired in the Newhouse-Robinson
Theorem, we construct a one-parameter family
$\{\widetilde{f}_t\}_{t\in I}$ perturbing in an adequate way the
diffeomorphism $f_0$. This perturbation is constructed so that there exists
a nested sequence of intervals of values of the parameter such that
in the $i+1$-interval there exists a sink $s_{i+1}$ and the $i$ sinks
constructed in the former intervals still persist. In the intersection of all
these intervals we obtain a parameter $t_\infty$ in which there
exist infinitely many sinks.




This construction is possible because $f_0$ has a homoclinic
tangency and perturbing $f_0$, a horseshoe is created. Newhouse
remarked the persistence of homoclinic tangencies of saddle points
of a horseshoe whose unstable and stable Cantor sets $K ^u$ and
$K^s$ along the line of tangencies have large thickness. Since near
a homoclinic tangency there exists a sink (see the Yorke-Alligood
theorem, (\cite{yorke})), it is possible to reason inductively in
order to construct the nested sequence of intervals.

Our purpose to prove Theorems \ref{teo}  and \ref{teoprevio} in this
paper, goes beyond the construction of Newhouse:   we shall  be
able, besides,  to perturb the primary family of diffeomorphisms in
the functional space $\mbox{Diff}^3(M)$, considering what we call
\lq \lq secondary diffeomorphisms",  along a properly defined
manifold ${\cal M}\subset \mbox{Diff}^3(M)$, in such a way that  the
infinite sinks, constructed for the diffeomorphism in the primary
family, persist simultaneously.

Taking a nearby family in an adequate infinite dimensional set of
$\mbox{Diff}^3(M)$, we will prove that the values of the parameter
where the tangencies and the sinks are produced, are near those of
the original family, and then the sinks continue, obtaining in this
way the manifold $\mathcal{M}$.



To prove part A of Theorem \ref{teoprevio}, with a suitable change of coordinates,
the diffeomorphisms of the family are near the  functions of the
classical quadratic family. For certain functions of the quadratic
family, the sink has eigenvalues as contractive as wanted. This
property is maintained for the diffeomorphisms of the original
family.

Part B will be proved  perturbing the diffeomorphism in such a way
that the sinks, which are far from a bifurcation after part A,
persist.

To prove part C we will show that any perturbation of the
diffeomorphism generically unfolding the quadratic tangencies allows
to persist only a finite number of the sinks, because the range of
the values of the parameter for which the sinks persist decreases
monotonically to 0.

The main tools that will allow us to make such proofs are
Propositions \ref{LemaBasico} and \ref{proposicion} of this paper.
This last guarantees  the existence of uniform sized manifolds of
codimension one in some infinite dimensional subset ${\cal N}_1
\subset \mbox{Diff}^3(M)$, along which all sinks persist
simultaneously.

As far as we  prove  our theorems, we construct the manifold ${\cal
M}$ having infinite dimension and also  infinite codimension. We do
not assert that the manifold $\cal M$ that we construct is maximal
verifying the conditions (a) and (b) of  the thesis of Theorem
\ref{teo}. Nevertheless,  if such a maximal manifold exists, it must
have at least codimension one, as a consequence of the part (C) of
the thesis in  Theorem \ref{teoprevio}, which we prove at the end of
the paper.

We also answer to other open question: Can the infinitely many
simultaneous sinks exhibited by a diffeomorphism $g_0$ constructed
as in Newhouse-Robinson Theorem simultaneously continue in an open
set? In fact, we prove that the answer is negative, provided that
$g_0$ is a diffeomorphism constructed  as in Theorem
\ref{teoprevio}.

\vspace{.4cm}

As a consequence of the proof of part (C) of Theorem \ref{teoprevio}, it is
immediate the following last result:

\begin{cor}
To continue infinitely many sinks \em (from those in
$\{s_i(f_\infty)\}$) \em of a diffeomorphism $f_\infty$ constructed
as in Theorem \ref{teoprevio} it is necessary to move along their
respective stable local leaves all the points in the unstable Cantor
set $ K^{u} $ of the line of tangencies $L(f_\infty)$ of $f_\infty$,
that are in the parabolic unstable arcs  of the accumulation points
of the sequence of sinks.
\end{cor}
This last result is the main reason why we restricted our
constructions (to prove the theorems of this paper) to an infinite
codimension manifold of diffeomorphisms, obtained from $f_\infty$
perturbing only inside $V$. In that way we can control easily the
unstable Cantor set $K^u$ while the stable Cantor set remains fixed.

\vspace{.3cm}

Finally, we pose the following open question. Let $f_\infty$ and
$\mathcal{M}$ verifying parts (A),  (B)1 and (B)2 of Theorem
\ref{teoprevio}. Has $\mathcal{M}$ necessarily infinite codimension?
%
%
%
%
%
%
%
%
%

\section {Persistence of tangencies.}

We recall the definition of line of tangencies (see Definition
\ref{df}) and stable and unstable Cantor sets (see Definition
\ref{Cantorsets}).

\begin{defi}
\em Given a Cantor set $K \subset \mathbb{R}$, the thickness  at
$u\in K$  in the boundary of a gap $U$ is defined as $\tau(K,u)=
\frac{l(C)}{l(U)}$ where $C$ is the bridge of $K$ at $u$ (see
\cite{palistakens}). The thickness of $K$, denoted $\tau (K)$ is the
infimum of the $\tau(K,u)$ over $u$.  \em
\end{defi}

\begin{defi} Two Cantor sets, $K_1$, $K_2\subset \mathbb{R}$ have large
thickness if $\tau(K_1)\tau(K_2)>1.$\em
\em
\end{defi}

\begin{defi} \label{espesurasgrandes} The horseshoe $\Lambda$ verifies
the large thickness condition if $\tau(K^s) \tau(K^u) >1$ where
$K^u$ and $K^s$ are defined in \ref{Cantorsets}. \em
\end{defi}

The importance of this definition resides in the following lemma:

\begin{lema}
Let $K_1$, $K_2 \subset \mathbb{R}$ two Cantor sets with large
thickness. Then, one of the following three alternatives occurs:
$K_1$ is contained in a gap of $K_2$; $K_2$ is contained in a gap of
$K_1$; $K_1 \bigcap K_2 \neq \emptyset$.
\end{lema}

For a proof, see \cite{palistakens}. We will apply the lemma to the
stable and unstable Cantor sets on the line of tangencies: the third
alternative assures the persistence of tangencies.

Now, we will define strongly dissipative horseshoes. Let us consider unstable and stable foliations $\mathcal{F}^u$ and $\mathcal{F}^s$ defined in a neighborhood $U$ of a horseshoe $\Lambda$. Let us take nonzero $C^1$ vector fields $X^u$ and $X^s$, tangent to the leaves of $\mathcal{F}^u$ and $\mathcal{F}^s$, and let us define the functions $\lambda\colon U \to \mathbb{R}$ and $\sigma\colon U \to \mathbb{R}$ as:
$$df(X^u(x))= \sigma(x) X^u(f(x))$$ $$df(X^s(x))= \lambda(x) X^s(f(x))$$
Redefining $X^u$ and $X^s$ if necessary, we obtain that $\sigma >1$ and $\lambda < 1$ for every point of $U$ if $U$ is a small enough neighborhood of $\Lambda$.
\begin{defi} \label{disipafuerte} \em
We say that a horseshoe is strongly dissipative if for every $x \in U$, $\lambda(x) \sigma^2(x) < 1$
\em
\end{defi}

\begin{teoremadeotro}
\label{newhouse} {\em \bf (Newhouse-Robinson)} Let $\widetilde{f}_t$
be a monoparametric family which generically unfolds a qua\-dratic
homoclinic tangency $q_0$ exhibited at $t=0$ of a fixed dissipative
saddle point  $P_0$ (i.e. the eigenvalues of $df_0(P_0)$ are
${\lambda}_0 < 1 < {\sigma}_0$ and ${\lambda}_0 {\sigma}_0 < 1$.)

Then, given $\varepsilon >0$, there exists an interval $I\subset (0,
\varepsilon)$ of values of the parameter and an open set $V$ such
that:
\begin{description}
\item[(i)] For every $ t \in I$ the diffeomorphism $\widetilde{f}_t$
exhibits a horseshoe $\Lambda$ which verifies the condition of large
thickness  as in Definition \em \ref{espesurasgrandes} \em and it is
strongly dissipative (i.e. $\lambda{(x)}  \sigma^2{(x)}<1, \; \;
\forall x \in \Lambda$).

\item[(ii)] For every value $t$ of the parameter in a dense set in $I$,
there exists a saddle point $P  \in \Lambda$ which exhibits an
homoclinic tangency $q \in  V$.

\item[(iii)] For every value $\tau$ of the parameter in a dense
set in $I$ there exists infinite simultaneous sinks in  $V$.
\end{description}
\end{teoremadeotro}

{\em Proof: } See \cite{newhouse} and \cite{robinson} and the lemma below.


The horseshoes created by unfolding tangencies have an important property:

 \begin{lema}  \label{disipacionfuerte}
The horseshoe created by the unfolding of a homoclinic quadratic
tangency of a dissipative saddle point can be taken strongly
dissipative just taking the number of iterates large enough.
 \end{lema}

{\em Proof: } It is a consequence of the scaling in Section 4,
Chapter III, of \cite{palistakens}. The horseshoe is
diffeomorphically conjugated to a horseshoe near a map of the
quadratic family which is infinite contractive ($\lambda =0$) along
its stable foliation. $\;\;\Box$

\vspace{.2cm}

To prove our first main result in Theorem \ref{teo} it is enough to
join the statements of Theorem \ref{newhouse} with the following:

\begin{teoremadeotro}
\label{teoespesura}

Let $M$ be  a  $C^\infty$ two dimensional compact connected
riemannian manifold. Let \em $f_1\in \mbox{Diff}^3(M)$ \em
exhibiting a strongly dissipative horseshoe $\Lambda \subset M$. Let
$P_1 \in \Lambda$ be a saddle periodic point  with a quadratic
homoclinic tangency at $q_1$.

Assume that the stable and unstable Cantor sets of $\Lambda$ along
the line of tangencies $L$ in a neighborhood $V$ of $q_1$ verify the
condition of large thickness defined in \em \ref{espesurasgrandes}.
\em

Then, given an arbitrarily small neighborhood ${\cal N}$ of $f_1$ in
\em $\mbox{Diff}^3(M)$, \em there exists a $C^1$ arc-connected
infinite-dimensional
manifold $\mathcal{M} \in {\cal N}$   such that:
\begin{description}
\item[(a)] Every $g\in
\mathcal{M}$ exhibits infinitely many simultaneous sinks $\{s_i(g)\}_{i \in \mathbb{N}} \in V$

\item[(b)] Each sink $s_i(g) $ is the   continuation of the respective sink $s_i(g_0)$
for any pair of diffeomorphisms $g_0, g \in \mathcal {M}$.
\end{description}

\end{teoremadeotro}

The proof of this Theorem is in Section \ref{seccionprueba1}.


\section {Local coordinates.} \label{seccionCoordenadas}

\vskip.5cm We  remark some known facts on the existence of a regular
coordinate system in a neighborhood $U$ of the horseshoe $\Lambda$
that trivialize its local stable and unstable  foliations.

{\obs  \label{coord} Regularity of the local invariant foliations.} \em

\em Given the horseshoe $\Lambda$ in a two dimensional manifold $M$,
and given a sufficiently small neighborhood ${U} \subset M$ of
$\Lambda$, there exist  the stable local foliation $\mathcal{F}^s$
and the unstable local foliation $\mathcal{F}^u$ that are invariant
while their  iterates remain in ${U}$ (see Appendix 1 in
\cite{palistakens}; see also \cite{demelo}). Moreover, if $f \in
\mbox{Diff}^2(M)$, then both invariant local foliations are of $C^{1
}$-class (see \cite {palistakens} Chapter II Section 3 and also
Appendix 1.) Then the stable leaves are $C^2$ and the tangent space
of the stable leaf through a point $P \in U$ depends $C^{1}$ on the
point $P$. In particular the concavity of each leaf depends
continuously on the point $P$.

Besides, if $f\in \mbox{Diff}^3(M)$ and $\lambda \sigma^2 <1$ then
the local stable foliation is of $C^3$-class while the unstable
foliation is not necessarily more then $C ^{1+\varepsilon}$. In
fact,  the $C^3$ differentiability of the stable foliation follows
after its $r-$normality (see \cite{hirshpughshub}): arguing as in
\cite{palistakens} Appendix 1, and working in the space $L(M)=
\left\{ (x,L):\; x \in M \right.$ and $L$ is a 1-dimensional linear
subspace of $\left. T_xM \right\}$, it follows that the local stable
foliation is $C^r$ with $r$ such that $\fra{\sigma \lambda^{-1} }{
\sigma^r} > 1$ for all $x$, or, equivalently, $\displaystyle r < 1 +
\frac{- \log \lambda}{\log \sigma}$. The $C^3$ regularity of
$\mathcal{F}^s$ follows recalling that we assumed that $\lambda
\sigma^2 < 1$.

We will mainly work with $C^3$ diffeomorphisms, so the stable
foliation will be $C^3$, and the unstable foliation, will be
$C^{1+\varepsilon}$.

\begin{obs}
{Local coordinate system in the neighborhood $U$ of the horseshoe. } \label{coord2}

\em As a consequence of Remark \ref{coord}, if $f \in
\mbox{Diff}^3(M)$ and $\lambda \sigma^2 <1$, we can take $C^{1}$
local coordinates $(x,y)$ of the two-dimensional manifold in the
neighborhood ${U}$ containing the horseshoe $\Lambda$, such that the
local stable leaves of $\Lambda$ are horizontal lines $y=
\mbox{constant}$ and its local unstable leaves are vertical lines
$x= \mbox{constant}$.

We get  $f(x,y) = (\xi (x), \eta (y))$ with $\xi $ of $C^1$ class and $\eta$ of $C^3$ class.

Also, given any $C^3$ map $H\colon U \mapsto U$ its computation in
the local coordinates $H(x,y)= (H_1, H_2)$ will be of $C^{1}-$class
and besides the second  and third order partial derivatives of $H_2$
respect to $y$ exist and are continuous.

Such a regular coordinate system in ${U}$ exists for any  map $g$ in
a neighborhood ${{\cal N}} \subset \mbox{Diff}^3(M)$ of the given
map $f_1$: in fact, the hypothesis of existence of the hyperbolic
horseshoe $\Lambda (g) $ verifying $\lambda \sigma^2
<1$, is persistent under small perturbations of $f_1$.

Besides, the local unstable and stable foliations and their tangent
spaces depend continuously on $g \in {\mathcal N}$. Therefore the
local coordinate system in ${U}$ chosen as above for each $g \in
{\mathcal N} \subset \mbox{Diff}^3(M)$, depends continuously on $g$.
\end{obs}


\begin{parag}
 \label{paragrafo}  The local coordinates computation of the map.
 \end{parag}
Take $f_1 \in \mbox{Diff}^3(M)$ exhibiting a horseshoe $\Lambda$ as
in the hypothesis of Theorem \ref{teoespesura} and a small neighborhood
${U}$ of $\Lambda$ in $M$.

We will perturb $f_1$ with diffeomorphisms $\xi$, i.e. $f=\xi \circ
f_1$ such that $\xi =$Id in a neighborhood of $\Lambda$, so that the
horseshoe remains the same. Later, we will choose an adequate
sequence of periodic saddle points $P_i \in \Lambda$ and use the
following notation:

\begin{obs}
\label{secuenciafsubi} Notation: \em

$P_i = P_i(f) \in \Lambda(f)$ is a saddle point. Let us denote
$\{f_i\}_{i \geq 1}$  a sequence of   diffeomorphisms $f_i \in
 \mbox{Diff}^3(M)$ along a one-parameter family from
$f_1$, such that $f_i$ exhibits a homoclinic quadratic tangency at
$q_i$ of the saddle $P_i$. Each of   the points $q_i$ shall be
chosen in a certain horizontal arc $A^s_i \subset \{y= y(q_i)\}\subset V$ of the
stable manifold of the saddle $P_i$, but not necessarily in the
local connected component $y= y(P_i)$ through $P_i$. The tangency
points $q_i$, for all $i \geq 1$,  shall be chosen in the line of
quadratic tangencies $L= L(f_i)$ contained in the
small open set ${V}\subset U $ defined in remark \ref{defincionV}.
\end{obs}

The horizontal arc $A^s_i \ni q$ is chosen small enough such that
$f^n_i(A^s_i)\subset U \, \forall n \geq 0$. Then we choose $n_i$
large enough so that $f^{n_i}_i(A^s_i) \subset \{y= y(P_i)\}$. Let
us take a height $h$ of a vertical segment $I_h$ such that
$f^j_i(A^s_i\times I_h) \subset U$ for $j=0$, $1,\ldots, n_i$. Now
we choose a region $D \subset f_i^{n_i}(A^s_i\times I_h)$ which
projects in a fundamental domain on $W_{\mbox{\footnotesize loc}}^u
(P_i)$. Finally, we take $V_i = f_i^{-n_i}(D)$, (let us observe that
$n_i$ can be taken as large as wanted) where we will rescale the
coordinates in the next Section.

Consider any $f_i $ as above. We will argue as in  \cite{palistakens} Chapter III, Section 4:

Taking $n_i $ large enough, and observing that in that case the
number of iterates near $P_i$ can be taken as large as wanted, we
have  that the length contraction $\lambda^{(n_i)}_i(A_i^s) $ of the
horizontal compact arc $A_i^s$ of stable manifold of $P_i$ in $U$,
when applied $f_i^{n_i}$, and the expansion $\sigma^{(n_i)}_i(A_i^u)$ of
a (vertical)  compact arc $A_i ^u$ of its unstable manifold in $U$,
verify  $\lambda^{(n_i)} _i(A^s_i) < \widetilde{\lambda}^{n_i} <1$
and
$\sigma^{(n_i)} _i (A^u_i)> \widetilde{\sigma}^{n_i} >1$ (see Remark \ref{defincionV}).

We have for all $(x,y) \in V_i$
\begin{equation} \label{valorespropios}
   \begin{array}{rcl}
    f_i(x,y)  & = & (\xi(x), \eta(y))
    \end{array} \end{equation}
 $$\mbox{dist}^s (f_i(P_i),(\xi(x), \eta(y)))= \lambda_i(x)
\mbox{dist}^s(P_i,(x,y))$$  
$$ 
\mbox{dist}^u(f_i(q_i),(\xi(x), \eta(y)))= \sigma_i(y)
\mbox{dist}^u(q_i,(x,y))
$$ where dist$^s$ and dist$^u$ can be taken as the
distances along compact arcs of stable and unstable manifolds of
$\Lambda$ as follows: dist$^u(q_i,(x,y)) = |y-y_{q_i}|$;
dist$^s(q_i,(x,y)) = |x-x_{q_i}|$ and dist$^s(P_i,q_i) $ can be
taken for instance as the length of the compact stable arc between
$P_i$ and the homoclinic tangency $q_i$ (we remark that this arc is
not necessarily contained in $U$, see Figure \ref{lastfigure}).

{\begin{figure}
[h]\psfrag{qi}{$q_i$}\psfrag{Pi}{$P_i$}\psfrag{q}{$q$}
\psfrag{U}{${U}$} \psfrag{s}{$\sigma$}\psfrag{l}{$\lambda$}
\psfrag{fn}{$f^n$}\psfrag{fN}{$f^N$}
\begin{center}\includegraphics[scale=.25]{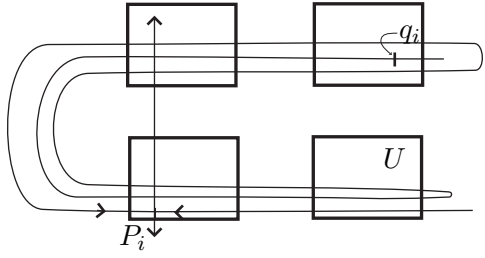}
\caption{\label{lastfigure} Stable distance.}\end{center}
\end{figure}}


If the point $(x,y)$ is such that $f_i^j (x,y) \in U$ for all $j=
0,1, \ldots, n_i$ and if $n_i$ is large enough:
    \begin{equation} \label{valorespropios2}
\sigma_{i} ^{(n_i)}(y)=  \left (
\prod_{j=0}^{n_i-1}\sigma_{i}(f_i^{j}(y)) \right ) > \widetilde{
\sigma} ^{n_i}, \; \; \; \; \; \lambda_{i} ^{(n_i)}(x)=  \left (
\prod_{j=0}^{n_i-1}{\lambda}_i(f_i^{j}(x)) \right ) <
\widetilde{\lambda}^{n_i}
   \end{equation}

In the last equations $f_i^{j}(y)$ denotes the ordinate of the point
$f_i ^{j}(x,y)$, which depends only on $y$. Similarly  $f_i^{j}(x)$
denotes its abscise, which depends only on $x$.

\vspace{.2cm}

\begin{parag}   The local coordinates computation near the quadratic homoclinic tangencies.
\label{paragrafo2}
\end{parag}
Let us study now the behavior of the map $f$ and its computation in
the linearizing coordinates near the homoclinic tangency. Let us
consider $f_i$ as in \ref{secuenciafsubi}, and the points
$r_i=(x(P_i),a_i)$ and $q_i=(b_i, c_i)\in V$ in the homoclinic orbit
of the saddle $P_i$, where the quadratic tangency is exhibited such
that $ q_i = f_i^{N_i} (r_i)$ for some integer number $N_i \geq 1$,
where $N_i$ is defined in the same way as $N_1$ in remark
\ref{defincionV}, but referred to the point $P_i$ .

Take $f \neq f_i$  in ${\cal N}$, and consider for such $f$ the
points $r_i(f)=(x(P_i(f)),a_i(f))$ and $q_i(f)=(b_i(f),c_i(f))$ as
in figure \ref{tg}, being the  \lq \lq remaining points" of  the
tangency that the map $f_i$ exhibited. These points $r_i(f)$ and
$q_i(f)$ are defined as follows:

First, we denote $r_i(f_i) = r_i, \; \; q_i(f_i)= q_i$. The point
$q_i$ belongs to a compact  arc $A_i^s$  of the stable manifold of
the saddle $P_i$, with equation $y= c_i$ when $f= f_i$. It is not
necessarily in the local stable manifold $y= y_{P_i}$ of the saddle
$P_i$.

Second, if $f \neq f_i $ in $ \mathcal{N}$ , then the homoclinic
tangency may disappear, but we still have the continuation
$A_i^s(f)$ of the compact stable arc $A_i^s$, with equation $y=
c_i(f)$,  and a new line of tangencies $L(f) \subset V$. We first
define the point $r_i(f)$ belonging to the connected component of
the local unstable leaf $x= x_{P_i}$ of $P_i$ in $U$, and being such
that $ f^{N_i}(r_i)\in K^u(f)$ is in the line of tangencies
$L{(f)}\subset {V}$.

Afterwards, we take the coordinates of the point $ f^{N_i}(r_i(f)) =
(b_i(f), \nu_i(f))=(b_i(f), \mu _i (f) + c_i(f))$. Its ordinate
$\nu_i(f)$ is the height of the parabolic arc in the compact piece
of unstable manifold of $P_i$ that made the tangency for the
diffeomorphism $f_i$. The height $\nu_i(f)$  is the sum of two
terms:    the \lq \lq relative height"   $\mu_i(f)$ respect to the
stable arc $A_i^s(f)$ with which that parabolic arc made the
tangency (i.e. $\mu_i(f_i) = 0$), and the ordinate $c_i(f)$ of the
arc $A_i^s(f)$. Finally we define the   point
$q_i(f)=(b_i(f),c_i(f)) \in A^s_i(f)$ as the projection of $
f^{N_i}(r_i(f))$ along the vertical direction on the stable leaf
$A_i^s(f)$, see Figure \ref{tg}.

 We remark that $r_i$, $q_i$, $a_i, b_i, c_i$ and $\mu_i$ depend continuously on $f \in \mathcal{N}$.

{\begin{figure}[h]\psfrag{P}{$P_i$}\psfrag{sigma}{$\sigma_{i}$} \psfrag{r}{$r_i=(x(P_i),a_i)$}
\psfrag{lambda}{$\lambda_{i}$}\psfrag{q}{$q_i=(b_i,c_i)$} \psfrag{m}{$\mu_i = \nu_i - c_i$}


\begin{center}\includegraphics[scale=.3]{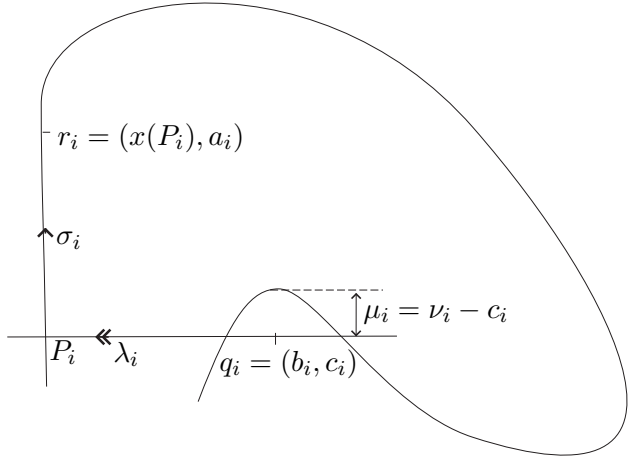}
\caption{\label{tg} Unfolding the tangency}\end{center}
\end{figure}}

We compute the equations of the transformation $f^{N_i}$   which
goes from a small neighborhood of $r_i = (x(P_i), a_i)$ to the
neighborhood $V$ of $q_i = (b_i, c_i)$ in the coordinates $(x,y)$.
We take $x^* = x - x(P_i), \; y^* = y- a_i$, and after
\cite{palistakens}:
$$
f^{N_i}\colon(x(P_i) + x^*,a_i+y^*) \mapsto(b_i,c_i)+ (H_1(\mu_i,x^*,y^*),H_2(\mu_i,x^*,y^*))
$$

Compute now the Taylor expansion of $H_1$ and $H_2$ in a
neighborhood of $(\mu_i, x^*, y^*)= (0,0,0)$, that is in a spacial
neighborhood of the point $r_i $ and a neighborhood of the
diffeomorphism $f_i$ in a one-parameter   family $\{f_{i, \mu_i}
\}\subset \mbox{Diff}^2(M)$ such that $f_{i,0} = f_i$ for  $\mu_i =
0$:

\begin{equation}
\label{alphabetagamma}
\begin{array}{rcl}
H_1(\mu,x^*,y^*) &=& \alpha_i y^* + \widehat{H}_1(\mu_i,x^*,y^*)\\
H_2(\mu,x^*,y^*)&=& \beta_i {y^*}^2 +\mu_i + \gamma_i x^* +\widehat{H}_2(\mu_i,x^*,y^*)\\
 & &\\
\mbox{where } \displaystyle \alpha_i= \frac{\partial H_1(0,0,0)}{\partial
y^*}, & &\displaystyle \beta_i = \frac{\partial^2 H_2(0,0,0)}{\partial
{y^*}^2}, \; \; \;\; \; \; \displaystyle\gamma_i= \frac{\partial
H_2(0,0,0)}{\partial x^*} \\
\end{array}
\end{equation}

\begin{obs}
\label{betacontinuo} \em Observe that the construction in the
subsection \ref{paragrafo2} is applicable for  $f \in$ Diff$^2(M)$.
Note that $\alpha _i, \beta_i, \gamma_i$ depend continuously on $f
\in $Diff$^2(M)$: in fact, the numbers $\alpha_i$ and $\gamma_i$ are
first order derivatives of the $C^{1}$- functions $H_1$ and $H_2$
which depend continuously on the given $f \in {\cal N}$. And
$\beta_i$ is a second order derivative along the stable foliation,
which is of $C^2$ class, and depends continuously on $f$, due to the
regularity of the chosen local coordinates and its continuous
dependence on $f$, as observed in Remark \ref{coord2}.

\em
\end{obs}

\begin{lema}
\label{alphabetanonulos}

$\alpha_i (f) \neq 0$, $\beta_i (f) \neq 0$ and $\gamma_i(f) \neq 0$
for all \em $f $ near enough $f_1$ in the $C^2$ topology.
\end{lema}

{\em Proof: } The result is due to the quadratic hypothesis. In
fact, $\partial H_2(0,0,0)/\partial y^* =0$ due to the tangency  at
the point $f^{N_i}(r_i)$. That is why we have chosen $r_i$ such that
$f^{N_i}(r_i)$ belongs to the line of tangencies $L_f \subset V$. As
$f_i$ is a diffeomorphism, the derivative $ Df_i(0,0,0) \neq 0$ and
so its determinant is not null, i.e. $\alpha _i \gamma_i \neq 0$.

Besides,  the tangency is quadratic and therefore Definition
\ref{quadratic} holds. Consider now the change of coordinates from
the $C^2$-system given in Definition \ref{quadratic} to the
coordinates leading to Equations (\ref{alphabetagamma}). We get the
relation $\overline {\beta_i }\neq 0$ (in Definition
\ref{quadratic}) if and only if $\beta _i /\alpha_i^2 \neq 0$ in
Equations (\ref{alphabetagamma}). Then $\beta_i \neq 0$ as wanted.
$\;\;\Box$


\begin{lema}
\label{lemaBeta}
If $\mathcal{N}$ is small enough, then for all $i $ there exists a real constant $K_i > 0$ such that for $f \in {\cal N}$ the coefficients $\alpha_i , \beta _i,
\gamma_i$ in equations \ref{alphabetagamma} verify:
$$ \frac {1} {K_i} \leq |\alpha _i|, \;|\beta_i |,\; |\gamma_i| \leq K_i$$
\end{lema}

{\em Proof:}


It is not restrictive to suppose a bounded small open set  ${\cal
N}\subset {\mbox{Diff}^3(M)}$, so for some $\varepsilon_0 >0$:
  \begin{equation}
  \label{a1}
  ||f-f_1||_{C^3}\leq \varepsilon_0\;\;\;\forall \, f \in {\cal {N}}  \end{equation}

We will  prove that there exists a positive lower bound  of
$|\beta_i|$ for all $f \in {\cal N}$. The proof of the existence of
the upper bound has a similar argument.

By contradiction, suppose that there exists a sequence of
diffeomorphisms $g_j \in {\cal N}\subset \mbox{Diff}^3 (M)$ such
that

 \begin{equation}
 \label{a2}
 |\beta_i (g_j)| \leq 1/j \; \; \forall \, j \geq 1
 \end{equation}

The sequence $g_j$ of diffeomorphisms in ${\cal N}$ is $C^3$-bounded
due to condition (\ref{a1}). By the Arzela-Ascoli Theorem there
exists a subsequence, which we still call $g_j$,  convergent in the
$C^2$ topology to a map $g_0 \in \mbox{Diff}^2(M)$. For this map
$g_0$ the number $\beta_i (g_0)$ in Equations (\ref{alphabetagamma})
is still defined and different from zero, due to Lemma
\ref{alphabetanonulos}.

As remarked in \ref{betacontinuo}, the real number $\beta_i (g)$
depends continuously on $g \in \mbox{Diff}^2(M)$. Therefore we get:
 $$ {\lim _{j \rightarrow \infty}}  {g_j}  = g_0 \; \Rightarrow \; |\beta_i(g_j)|
 \rightarrow |\beta_i(g_0)| \neq 0 $$

Therefore the sequence of real numbers $ |\beta_i(g_j)|$ is bounded
away from zero, contradicting the inequality (\ref{a2}). $\;\;\Box$


 \section{Approximation to the one-dimensional quadratic family.} \label{sectionqfamily}

 We continue arguing
 as in  \cite{palistakens} Chapter III, Section 4:

\vspace{.3cm}

Consider  $f \in {\cal N} \subset \mbox{Diff}^3(M)$ as in Section
\ref{seccionCoordenadas}, and for fixed $i\geq 1$ take the periodic saddle point $P_i
\in \Lambda$ and the coordinate system $(x,y)$ defined in Remark
\ref{coord2} in the neighborhood $U$ of $\Lambda$.

Take the point $r_i= (x(P_i), a_i)$ in the local unstable vertical
arc through $P_i$,  and the point $ q_i = (b_i, c_i)$ in the
horizontal leaf  $y= c_i$ contained in the global stable manifold of
$P_i$ in $U$, as defined in Section \ref{seccionCoordenadas} and
Figure \ref{tg}.

We recall equations (\ref{valorespropios}) and will consider a
change of coordinates in the small open rectangle  $V_{i} \subset U$
near $q_i$ defined in \ref{secuenciafsubi}.

The following change of variables, and also the reparametrization on
the value of $\mu_i$, are defined in \cite{palistakens} Chapter III,
Section 4, near a quadratic homoclinic tangency. We have made some
minor adaptation to our context, in which the coordinate system
$(x,y)$ in the neighborhood $U$ of the horseshoe $\Lambda$ is
independent of the saddle point $P_i \in \Lambda$ with which we
work. Therefore $P_i$ does not have necessarily coordinates $(0,0)$.
We write:


\begin{equation} \label{reparam}
\begin{array}{l}
\widehat{{\mu _i}}=\left( \mu_i [\sigma_{i} ^{({n_i})}(y)]^{2}
 + \mbox{dist}^s(q_i,P_i)
\gamma_i [\lambda_{i} ^{({n_i})}(x)]
 [\sigma_{i} ^{({n_i})}(y)]^{2}- (a_i - y(P_i))
[\sigma_{i}^{({n_i})}(y)]\right)\beta_i\\
   \widehat{x} = (x-b_i)
[\sigma_{i} ^{({n_i})}(y)]\beta_i \alpha_i^{-1}\\
    \widehat{y} = ((y - c_i)
    [\sigma_{i} ^{({n_i})}(y)]^{2}  - (a_i - y(P_i))
    [\sigma_{i} ^{({n_i})}(y)])\beta_i
  \end{array}
\end{equation}
where  the definition of the coefficients $\sigma_{i} ^{({n_i})}(y)$ and $ \lambda_{i} ^{({n_i})}(x)$ are in Equations  (\ref{valorespropios2}), and $x(Q)$ and $y(Q)$ denote respectively the abscissa and ordinate of $Q$. We recall that if ${n_i}$ is large enough then:
$$ \sigma_{i} ^{({n_i})}(y) >  \widetilde{\sigma} ^{n_i}, \; \; \; \; \; \lambda_{i} ^{({n_i})}(x)<   \widetilde{\lambda} ^{n_i} $$
where  $\widetilde{\lambda} <1$ and $\widetilde{\sigma} >1$ are the
exponential contractive and expansive rates of the hyperbolic set
$\Lambda$. The Inverse Function theorem allows us to assert that the
former equations define invertible $C^1$ change of coordinates. We
recall that at each point, $\lambda \sigma^2 <1$.

For later use we write the following equations, obtained from \ref{reparam}:

\begin{equation} \label{reparam3}
\begin{array}{l}
\mu_i=  \beta_i^{-1}{\widehat{\mu}_i} [\sigma_{i}
^{({n_i})}(y)]^{-2}  - \mbox{dist}^s(q_i,P_i)
\gamma_i [\lambda_{i} ^{({n_i})}(x)] + (a_i - y(P_i))
[\sigma_{i} ^{({n_i})}(y)]^{-1}  \\
x=b_i+ \alpha_i \beta_i^{-1}\widehat{x} [\sigma_{i} ^{({n_i})}(y)]^{-1}   \\
y= c_i + (a_i - y(P_i)) [\sigma_{i} ^{({n_i})}(y)]^{-1} +
\beta_i^{-1} \widehat{y} [\sigma_{i} ^{({n_i})}(y)]^{-2}
 \end{array}
\end{equation}

Given a point $(\widehat{x}, \widehat{y})$ in the new system of
coordinates, we apply $f^{n_i+N_i}$ (with $\mu$ constant), using
Equations \ref{valorespropios} for the first $n$ iterates of $f$,
and Equations \ref{alphabetagamma} for the last $N_i$ iterates. The
detailed computations are explicit   in \cite{palistakens} Chapter
III, Section 4. We get
$$ \left(
     \begin{array}{c}
       \widehat{x} \\
       \widehat{y} \\
     \end{array}
   \right) \stackrel{f^{n_i+N_i}}{\longrightarrow} \left(
   \begin{array}{c}
   F_1(\widehat{x},\widehat{y},\widehat{\mu}_i,{n_i}) \\
   F_2(\widehat{x},\widehat{y},\widehat{\mu}_i,{n_i}) \\
   \end{array}
   \right)
$$

The value of $\widehat{\mu}_i$ is obtained computing $(x,y)$ through
the last two equations of \ref{reparam3} and then substituting in
the first equation \ref{reparam}. We note from the first equation
\ref{reparam} that being $\mu_i$ constant, the value of
$\widehat{\mu}_i$ changes when applying $f^{n_i+N_i}$ because it
depends on $(x,y)$ which changes when applying the map.

For the next lemma, we consider in $D=[-1, 1]^3$  the 2-dimensional
manifold $\mathcal{S}$ of points
$(\widehat{x},\widehat{y},\widehat{\mu}_i)$ implicitly defined by
equations \ref{reparam} with a fixed value $\mu_i$. It can be
written as $\widehat{\mu}_i = g(\widehat{x},\widehat{y} )$. Let us
observe that for $n_i$ large enough, $\mathcal{S}$ approaches to a
horizontal surface:

\begin{lema} \label{nuevo2}
 $\frac{\partial \widehat{\mu}_i}{\partial \widehat{x}}$ and
$\frac{\partial \widehat{\mu}_i}{\partial \widehat{y}}$  converge
uniformly to 0 for $n_i\to \infty$ and $(
\widehat{x},\widehat{y},\widehat{\mu}_i)\in [-1,1]^3$.
\end{lema}

{\em Proof}
$$\left|\frac{\partial \widehat{\mu}_i}{\partial \widehat{y}}\right| =
\left| \frac{\partial \widehat{\mu}_i}{\partial x}\cdot \frac{\partial x}
{\partial \widehat{y}} + \frac{\partial \widehat{\mu}_i}{\partial y}\cdot
\frac{\partial y}{\partial \widehat{y}}\right| $$

Computing:

$$\frac{\partial y}{\partial \widehat{y}}= \frac{\beta_i^{-1}
[\sigma_i^{(n_i)}]^{-2}}{1+\left(\frac{a_i-y(P_i)}{[\sigma_i^{(n_i)}]^{2}}+
\frac{2\beta_i^{-1}\widehat{y}}{[\sigma_i^{(n_i)}]^{3}}\right)\sum_j\sigma^{(n_i-1)}_{i,j}
\sigma'_{i,j}}$$ where $\sigma_{i,j}^{(n_i-1)}$ is a notation for
the product $\sigma_i^{(n_i)}$ (see equation
(\ref{valorespropios2})) where we take out the $j$-th factor, and
$\sigma_{i,j}' $ is the notation for the derivative of the omitted
factor. It follows that there exists $k_i$ such that if $n_i$ is
large:
$$\left|\frac{\partial y}{\partial \widehat{y}}\right|\leq  k_i[\sigma_i^{(n_i)}]^{-2}$$

Similarly, 

$$\frac{\partial x}{\partial \widehat{y}} = -\,\frac{\alpha_i \beta_i^{-1}\widehat{x}}
{[\sigma_i^{(n_i)}]^{2}}\sum_j\sigma^{(n_i-1)}_{i,j}
\sigma'_{i,j}\frac{\partial y}{\partial \widehat{y}}$$

$$\frac{\partial \widehat{\mu}}{\partial x}= \mbox{dist}^s(q_i,P_i) \gamma_i \left(\sum_j\lambda_{i,j} ^{({n_i}-1)}\lambda'{i,j}\right)  [\sigma_{i} ^{({n_i})}(y)]^{2}\beta_i
$$

$$\frac{\partial \widehat{\mu}}{\partial y}=\left(2\left( \mu_i   + \mbox{dist}^s(q_i,P_i) \gamma_i [\lambda_{i} ^{({n_i})}(x)] \right) [\sigma_{i} ^{({n_i})}(y)]- (a_i - y(P_i))\right)\beta_i\sum_j\sigma^{(n_i-1)}_{i,j} \sigma'_{i,j}$$

Moreover, from the definition of strongly dissipative horseshoe and
the first equality of (\ref{reparam3}), it follows that there exists
$k_i$ large enough such that $|\mu_i | \leq k_i
(\sigma^{(n_i)}_{i})^{-2}$.  Therefore, increasing $k_i$ if
necessary,

$$\left| \frac{\partial \widehat{\mu}}{\partial y}\right|\leq k_i n_i [\sigma_{i} ^{({n_i})}(y)]$$

We will take numbers $\lambda^*$, $\lambda^+$, $\mu^*$, $\mu^+$ with $\lambda^* <\lambda(x)< \lambda^+ < 1<\sigma^*<\sigma(y)< \sigma^+$ $\forall \, (x,y) \in V_i$ such that $\lambda^+ \sigma^+< 1$.

$$\left|\frac{\partial \widehat{\mu}_i}{\partial \widehat{y}}\right|\leq k_i(n_i^2(\lambda^+\sigma^+)^{n_i}(\sigma^*)^{-2n_i}+n_i {\sigma^*}^{-n_i}) \longrightarrow_{n_i \to \infty} 0$$
%
uniformly in  $(\widehat{x},\widehat{y},\widehat{\mu}_i) \in
[-1,1]^3$ as wanted. Analogously it is proved for $\displaystyle
\frac{\partial \widehat{\mu}_i}{\partial \widehat{x}}$. $\;\;\Box$

We  conclude that taking $n_i\to \infty$, $f^{n_i + N_i}|_{V_{n_i}}$
converges in the $C^1$ topology, uniformly to the asymptotic map:

%
%
%
%
%
\begin{equation}\label{qfamily}
\left( \begin{array}{c}\widehat{x} \\ \widehat{y}\end{array} \right) \mapsto \left( \begin{array}
{c} \widehat{y} \\ \widehat{y}^2 + \widehat{\mu}\end{array} \right)
 \end{equation}

\begin{obs}
\label{pozo0} \em Note that the family defined by  Equation
\ref{qfamily} is the one-dimensional quadratic family with parameter
$\widehat{\mu}$. It is standard to verify that this quadratic family
exhibits a fixed point which is a sink for the parameter values
$\widehat{\mu} \in \left(-\,\frac34, \frac14\right)$
and that its basin of attraction includes all points $(\widehat{x},\widehat{y})$ with $ \widehat{y}$ in the interval $ (-1/4,1/4)$. 

Even more, if $\widehat \mu < -3/4$ of  if $\widehat \mu >1/4$, the
fixed point in the one-dimensional quadratic map  does not exist or it
is not a sink.

Observe that for the one-dimensional  quadratic family, the sink has
two eigenvalues: one is always zero, along the horizontal lines $y
=$constant, because it has infinite contraction transforming the
horizontal line onto one single point. The other eigenvalue  is the
slope at the sink of the parabola $\widehat y \mapsto \widehat y^2 +
\widehat \mu$.

We note that for $\widehat \mu = -3/4$  the sink has an  eigenvalue
equal to $-1$ and the quadratic unidimensional family exhibits there
a period doubling bifurcation. On the other hand, if $\widehat \mu =
1/4$ the sink has eigenvalue equal to 1, and the family has a saddle
node bifurcation. For $\widehat \mu \in (-3/4, 1/4)$ the slope of
the parabola at the sink, (being less than 1 in absolute value), is
continuous and monotone with $\widehat \mu$. Therefore, given any
$0<\rho<1$, there exist numbers $-3/4 < 2k ^-(\rho)<0< 2k^+(\rho) <
1/4$ such that if $\widehat  \mu\in (2k ^-(\rho), 2k^+(\rho))$ then
the  sink has both eigenvalues smaller than $\rho$ in absolute value.
\end{obs}

After the changes of coordinates and the reparametrization given in
Equations (\ref{reparam}),  $f^{n_i + N_i} $ converges uniformly to
the quadratic family when $n_i \rightarrow + \infty$. The speed of
convergence depends on the values of the hyperbolic expansive rates
$\sigma(f), \lambda(f)$ in the horseshoe exhibited by $f$ and also
of the values of $ \alpha_i(f), \beta_i(f), \gamma_i(f), ||\widehat
H_1 (f)||_{C^0}, ||\widehat H_2(f)||_{C^0}$, defined in equations
(\ref{alphabetagamma}). Due to Lemma \ref{lemaBeta},  these are
uniformly bounded for all $f \in {\cal
N}$.  

{\begin{figure} [h]\psfrag{Pi}{$P_i$}\psfrag{yqi}{$y_{i}$}
\psfrag{yri}{$y_{i+1}$} \psfrag{Qi}{$P_{i+1}$}\psfrag{bi}{$c_{i}$}
\psfrag{mui}{$\nu_{i}$}\psfrag{muj}{$\nu_{i+1 }$}\psfrag{muv}{$\mu_{i }$}
\psfrag{po}{$P_1$}\psfrag{ai}{$a_{i}$}
\psfrag{yo}{$y_{0,t}$}\psfrag{muo}{$t=\mu_{1}$}
\begin{center}\includegraphics [scale=.35,angle=270]{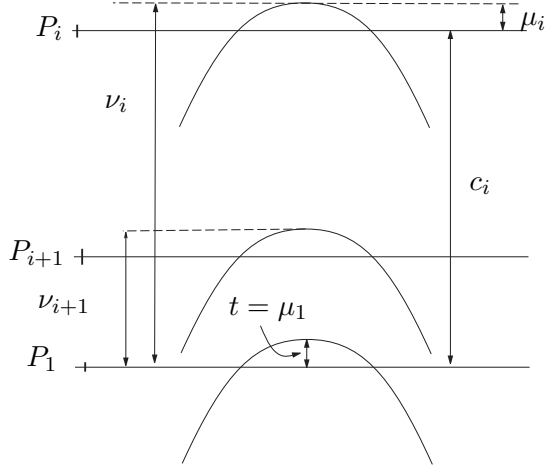}
\caption{\label{gr10} Construction of sinks.}\end{center}
\end{figure}}

We recall the notation of subsection \ref{paragrafo2} and Figures
\ref{tg} and \ref{gr10}. We choose a large enough natural number
$n_i$ and the small rectangle $V_{i} \subset V$, defined in Remark
\ref{secuenciafsubi}, to apply Equations (\ref{reparam}) which lead
asymptotically to the one-dimensional quadratic family, and apply
the results in Remark (\ref{pozo0}). We now resume in Proposition
\ref{LemaBasico} all the conclusions obtained in this section.

 \begin{prop}
\label{LemaBasico}
\begin{description}
\item[(A)] If for some real number $\widehat \mu \in (-3/8, \; 1/8) $ and for some $n_i$ large enough is verified
 \begin{equation}
 \begin{array}{rcl}
 \label{lemabasicoigualdad}
 \nu_i(f)& =&c_i(f)  + \mu_i(f)=\\
 &=& c_i(f) + (a_i(f)- y(P_i)) \cdot  [\sigma_{i}(f)^{(n_i)}(\nu_i(f)]^{-1} +\\ && + \widehat \mu  \cdot {\displaystyle \frac{[\sigma_{i}(f)^{(n_i)}(\nu_i(f))]^{-2}}{\beta_i(f)}}
  - \mbox{\em dist}^s(q_i,P_i) \cdot  \gamma_i(f) \cdot  [\lambda_{i}(f)^{(n_i)}(b_i(f))]\\
 \end{array}
 \end{equation}
  then $f$ exhibits a sink $s_i = s_i(f)$ in the given open set $V$.


\item[(B)] If $f, g \in {\cal N}$ are arc-connected in ${\cal N}$ by a one-parameter family $\{\widetilde{f}_t\}_t $ that verifies  the equality \em (\ref{lemabasicoigualdad}) \em for some  $C^1$ real function $$ \widehat \mu (\widetilde{f}_t) \in (-3/8, \; 1/8) \; \; \forall t$$ then the sink $s_i(g)$ is the   continuation of the sink $s_i(f)$.

\item[(C)] \label{D} If $f$ and $g$ verify equality \em(\ref{lemabasicoigualdad}) \em for some  $\widehat \mu (f) \in (-3/8, \; 1/8)$  and some
$\widehat \mu (g)< -1$ or $\widehat \mu (g) >1$,  then there does
not exist the   continuation of the sink $s_i(f)$ for such $g$.

\item[(D)] Given $0 <\rho <1$ there exist constants $-3/8< k^-(\rho) < 0 < k^+(\rho) < 1/8$ such that if for some real number $\widehat \mu \in (k^-(\rho), k^+(\rho)) $ and for some $n_i$ large enough is is verified equation (\ref{lemabasicoigualdad})
then $f$ exhibits a sink $s_i = s_i(f)$ in the open set $V$
with both eigenvalues smaller than $\rho$.
\end{description}
\end{prop}

{\em Proof: } Recall   that $\nu_i(f) = c_i(f)  + \mu_i(f)$ and take
into account (\ref{reparam3}). The real number $\widehat{\mu}$ in
equation  (\ref{lemabasicoigualdad}) is the parameter
$\widehat{\mu}_i$ in the first equality of (\ref{reparam3}) for
$x=b_i(f)$ and $y=\nu_i(f)$. As $\widehat{\mu} \in \left(-3/8,
1/8\right)$ and $n_i$ is large enough the reparametrized map $f^{n_i
+ N_i}$ is uniformly near the quadratic family and due to the remark
\ref{pozo0}, it exhibits a sink.

Part B is a consequence of the Implicit Function Theorem applied to (\ref{reparam3}).

Parts C and D follow after the remark  \ref{pozo0}. $\;\; \Box$
%
%
%
%


\section{Uniform sized   continuation of the sinks.} \label{seccionN}

Given $f_1 \in \mbox{Diff}^3(M)$ verifying the hypothesis of Theorem
\ref{teoespesura} let us consider the neighborhood ${\cal N} \subset
{\mbox{Diff}^3}(M)$ of $f_1$ as in Section \ref{seccionCoordenadas}.

Take a small neighborhood $U$ of the horseshoe $\Lambda$ and the
coordinate system as in Section \ref{seccionCoordenadas}. In this
Section  we shall assume the hypothesis of the strong dissipative
horseshoe $\lambda \sigma ^2 < 1$ so the local coordinate system
that trivializes the foliation of $\Lambda$ is of $C^3$ class.

Consider the homoclinic tangency point $q_1 = (b_1,0) \in U$ of the
saddle $P_1 = (0,0) \in \Lambda$ as in the hypothesis of Theorem
\ref{teoespesura},  and the line of tangencies $L_1 = L(f_1) \ni
q_1$ in $V$.

Let us define the following $C^{1}$ manifold ${\cal N}_1 \subset
{\cal N}\subset {\mbox{Diff}^3(M)}$,    which has infinite dimension
and codimension in $\mbox{Diff}^3(M)$, contains $f_1$, and will be
considered our universe where $f_1$ shall be perturbed.

For any $\delta >0$ small enough, (to be fixed later) we define:
\begin{equation}
\label{ambienteRestringido}
{\cal N}_1 = \left \{ f \in {\cal N}: f= \xi \circ f_1    \right \}\\
\end{equation}

where $ \xi \in \mbox{Diff}^3(M) , \; \; ||\xi - id||_{C^3} < \delta $   and besides:
\begin{equation}
\label{ambienteRestringido2}
   \begin{array}
 {l}
 \xi (p) = p \; \forall p \not \in V   \\
\exists k= k(\xi)\in \mathbb{R}   \mbox{ such that }\forall \, (x,y)
\in L_1: \\ \pi_2( \xi (x,y)) = y + k,  \; \; \;
(D{(\pi_2 \circ \xi)}(x,y)) \cdot (1,0) = 0 \; 
    \\
 \end{array}
\end{equation}
Here we denote  $\pi_2$ to the horizontal projection $\pi_2(a,b)=b$,
which is of $C^3$ class, due to the choice of the coordinate system
in Remark \ref{coord2}, under the assumption of the strong
dissipative hypothesis $\lambda \sigma ^2 < 1$.

We observe that $\xi $ is isotopic with the identity map, so
${\cal N}_1 $ is an arc connected manifold.


 \begin{lema} \label{variedadCod1}
For $K $ small enough  the following set
 $${\cal M}_K = \{f= \xi \circ f_1 \in {\cal N}_1: k(\xi)= K \}$$
   is a $C^1$ submanifold of codimension one in ${\cal N}_1$.

\end{lema}
{\em Proof:}

The map $k: \; \xi  \mapsto k(\xi)$, defined for all $\xi \in
\mbox{Diff}^3(M)$ which verify the conditions
(\ref{ambienteRestringido2}), is the second coordinate of the vector
obtained by the evaluation of $\xi - Id$ at  $q_1= (b_1,0) \in L_1=
L(f_1)$. Therefore, the map $k$ is a $C^1$ real function defined in
the set of diffeomorphims $\xi$ verifying conditions
(\ref{ambienteRestringido2}). Besides, its Fr{\'e}chet derivative
respect to $\xi$ is the evaluation $\xi \mapsto \xi (q_1)$ which is
a not null linear transformation on $\xi$ in the tangent space of
${\cal N}_1$. Therefore the value $K$ is a regular value of the real
function $k$, and so the equation $k(\xi) = K$ for $\xi$ such that $
\xi \circ f_1 \in {\cal N}_1$ defines a $C^1$ submanifold ${\cal
M}_K \subset {\cal N}$ of codimension one in ${\cal N}_1$, as
wanted. $\;\;\Box$

\vspace{.3cm}

The conditions (\ref{ambienteRestringido}) and
(\ref{ambienteRestringido2})  mean that we are perturbing $f_1$
only in the neighborhood $V$ near the line of tangencies $L_1 =
L(f_1)$ in such a way that we apply a vertical translation of
amplitude $k(\xi)$ to $L_1$ to obtain the new line of tangencies
$L(f)$ for $f = \xi \circ f_1$, and a horizontal deformation.

In particular we neither perturb the horseshoe  $\Lambda$, nor the
diffeomorphism in a neighborhood of $\Lambda$. Therefore, the local
stable and unstable manifolds of $\Lambda$ are the same, and the
system of local coordinates in $U$, as defined in Section
\ref{seccionCoordenadas}, does not change when perturbing $f$.

For each $i\geq 1 $ we choose any sequence of periodic saddle points
$P_i \in \Lambda$ as in the subsection \ref{paragrafo}. Let us
suppose a one-parameter family of diffeomorphisms in ${\cal N}_1$
having a  sequence of  diffeomorphism $f_i \in {\cal N}_1$ which
exhibits a homoclinic tangency at $q_i $ of the saddle $P_i$. We use
the notation of subsection \ref{secuenciafsubi}.

We are working along the restricted  space ${\cal N}_1$ of
diffeomorphisms that coincide with  $f_1$ in a neighborhood of the
horseshoe $\Lambda$. Therefore,  the values $\lambda_{i} (x),
\; \sigma_{i} (y)$ in Equation (\ref{valorespropios}) and of
$\lambda_{i}^{(n)}(x)$ and $\sigma_{i}^{(n)}(y)$ in Equations
(\ref{valorespropios2}) and (\ref{reparam}), are the same for all $f
\in {\cal N}_1$.

When passing from $f_1$ to $f \in {\cal N}_1$, the horizontal local
stable foliation remains fixed, and we apply a  transformation $\xi$
preserving the horizontal direction in the points of the line of
tangencies $L_1$ to obtain $L (f)$. Then, the point $q_i = q_i(f)
=(b_i(f),c_i)$ remains in the same horizontal line ($c_i$ is fixed)
and the point $f^{N_i}(r_i) = (b_i(f), c_i + \mu_i(f)) \in L(f)$
moves from $f^{N_i}_1(r_i) = (b_i(f_1), c_i + \mu_i(f_1)) \in
L(f_1)$ a vertical distance $k= k(\xi) = k(f \circ f_1^{-1})$, and
slides horizontally preserving its quality of being a point in the
line of tangencies. Therefore, the numbers $a_i,  \; c_i$ and  the
point $r_i$, defined in subsection \ref{paragrafo2} and Figure
\ref{tg}, remain the same for all $f \in {\cal N}_1$ and
$\mu_i(f)-\mu_i(f_1)=k(\xi) = k(f\circ f_1^{-1})$. In particular if
$f=f_i=\xi_i\circ f_1$ such that $\mu_i(f_i)=0$, we obtain
$-\mu_i(f_1)=k(\xi_i)$ and therefore $\mu_i(f) + k(\xi_i) = k(\xi)
\forall \, f \in \mathcal{N}_1 \forall \, i \geq 1$.

We conclude the following:
\begin{obs}
\label{constantesenespaciofuncional}
{\em
The points $P_i$ and $r_i$,  the real numbers $a_i $ and $ c_i$, and the functions $ \sigma_{i}^{(n) }$ and $
\lambda_{i}^{(n)}$ do not depend on the diffeomorphism $f \in {\cal N}_1$.}
\end{obs}

\begin{obs} {\em
$ \mu_1(f) = k(\xi)= k(f \circ f_1^{-1}) \forall \, f \in
\mathcal{N}_1$, and for $f_i = \xi_i \circ f_1 \in {\cal N}_1$ such
that $\mu_i (f_i)=0$ we obtain \begin{equation} \label{muconk}
\begin{array} {rcl}
k(\xi) &=&    k(\xi_i) + \mu_i(f) \forall \, f \in \mathcal{N}_1  \\
\nu_i(f) & = & c_i + \mu_i(f) =
c_i + k(f \circ f_1^{-1}) - k(f_i \circ f_1^{-1})  \forall \, f \in {\cal N}_1 \\
\end{array}
\end{equation}}
\end{obs}
\vspace{.3cm}

We recall that (\ref{ambienteRestringido2}) assumes $||\xi - id ||_{C^3} < \delta $, so we obtain the following:


\begin{lema}
\label{lemaNuevo}  Given $\varepsilon >0$ there exists $\delta >0$
such that if the manifold ${\cal N}_1 = {\cal N}_1 (\delta )$ is
constructed fulfilling Equations (\ref{ambienteRestringido}) and
(\ref{ambienteRestringido2}), then the number $ \beta _i(f)$ defined
by the equations (\ref{alphabetagamma}), verifies the following
inequalities  for all $f \in {\cal N}_1$ and all $i \geq 1$:

$$(1-\varepsilon ) |\beta_i(f_1)| < |\beta _i(f)| < (1+\varepsilon ) |\beta_i(f_1)|$$

\end{lema}

{\em Proof: } We have $f= \xi \circ f_1$, with $\xi (p) = p \; \;
\forall \, p \not \in V$, and $V$ the neighborhood of the line of
tangencies defined in Remark \ref{defincionV}.

Recall that for all $i$ the point $r_i$, and its first $N_i-1$
forward iterates, do not lay in $V$, and $f_1^{N_i}(r_i) \in V$. And
this also holds for all the points $p$ in a small open neighborhood
of $r_i$. As $f(p)= f_1(p) \; \forall \, p \not \in V$, we deduce
$f^{N_i} (p) = f \circ f^{N_i -1}(p) = \xi \circ f_1^{N_i}(p)$ for
all the points $p$ in a small open neighborhood of $r_i$.

We recall the definition in equalities (\ref{alphabetagamma}): the
number $\beta_i(f)$ is the second order partial derivative respect
to $y$ of the $C^3$ transformation $\pi_2 \circ f^{N_i} = \pi_2
\circ \xi \circ f_1^{N_i}$ at the point $r_i$.

We will work in a new $C^2$ system of coordinates (to be able to
apply the chain rule), such that the lines $y$ constant coincide
with the stable foliation in $U$ and we take a non invariant
foliation as $x$ constant. In these new coordinates the value of
$\beta_i(f)$ is the same as in the former system. The first
derivative respect to $x$ of $\pi_2 \circ \xi$ at $f_1^{N_1} (r_i)
\in L(f)$ is null due to our assumption that $D(\pi_2\circ \xi)(1,0)
= 0$ in the line of tangencies $L(f)$. On the other hand, its  first
derivative respect to $y$ is in $(1-\delta, 1 + \delta)$ due to $||
\xi - id ||_{C^3} < \delta$.

Now, denoting $(u,v)= f_1^{N_i}(x,y)$:

$$\frac{\partial(\pi_2  f^{N_i})}{\partial y}= \frac{\partial
(\pi_2 \xi)}{\partial u} \frac{\partial u}{\partial y} +
\frac{\partial (\pi_2 \xi)}{\partial v} \frac{\partial v}{\partial
y} $$ and then,  (we omit the points at which
we evaluate the partial derivatives) using that $\frac{\partial
(\pi_2 \xi)}{\partial u}= 0$, $\frac{\partial v}{\partial y}=0$ it follows:

$$\beta_i(f) =
\frac{\partial^2(\pi_2  f^{N_i})}{\partial y^2}=  \frac{\partial^2 (\pi_2\xi)}{\partial u^2}\left(\frac{\partial u}{\partial y}\right)^2+\frac{\partial
(\pi_2 \xi)}{\partial v} \frac{\partial^2v}{\partial y^2} $$

The  first term is bounded by $\delta \alpha_i^2(f_1)$, which can be
taken smaller than $\beta_i (f_1) \varepsilon/2$  taking $\delta<
\inf_i \{\frac{\varepsilon \beta_i}{2\alpha_i^2}\}$. We note that
$\frac{ \beta_i}{\alpha_i^2}$ is the concavity of the unstable
parabolic  arcs which are uniformly bounded away from 0 due  to the
quadratic hypothesis. The second term belongs to $((1-\delta)
\beta_i(f_1),(1+\delta) \beta_i(f_1))$. Then, taking $\delta <
\varepsilon/2 $ , $\beta_i(f)$ belongs to the interval
$((1-\varepsilon) \beta_i(f_1),(1+\varepsilon) \beta_i(f_1))$.
$\;\;\Box$

 \vspace{.2cm}



\begin{lema}
\label{Lemabasico2} For each $i \geq 1$ there exist   real constants
$m_i$, $K_i >0$ and $\nu_i^{(0)}$ (that are independent of $f\in
{\cal N}_1$) such that if $n_i > m_i$ then the implicit function
$\nu _i(f) (\widehat \mu) = G_f(\widehat \mu )$  defined by the
equation \em (\ref{lemabasicoigualdad}) \em  in Lemma \em
\ref{LemaBasico}, \em verifies:

$$|G_f(0) - \nu_i^{(0)}| \leq K_i (\max_{(x,y) \in \overline{V}_i}\{\lambda (x)\})^{n_i}$$
$$\frac{(\max_{(x,y) \in \overline{V}_i} \{\sigma  (y)\} )^{-2n_i}}{(1+\varepsilon)|\beta_i(f_1)| }
<|G_f' (\widehat \mu)| \approx
\frac{[\sigma_{i}^{(n_i)}(\nu_i)]^{-2}}{|\beta_i(f)|}<
\frac{\widetilde{\sigma} ^{-2n_i}}{(1-\varepsilon)|\beta_i(f_1)| }
\; \; \forall \, \widehat \mu \in (-2, 2)$$

Moreover the constant number $\nu_i^{(0)}$ is the sum of the
following two  terms, each independent of $f \in{\cal N}_1$:
\begin{equation}
\label{lemabasico2igualdad} \nu_i^{(0)}=c_i+ (a_i- y(P_i)) \cdot
[\sigma_{i}^{(n_i)}(\nu_i^{(0)})]^{-1}
\end{equation}

\end{lema}

\begin{figure}  [h]\psfrag{aim1}{$c_{i}$}\psfrag{nu}{$G_f(0)$}\psfrag{nu2}{$\nu_i^{(0)}$}\psfrag{nma}{$\nu_i^+$}
\psfrag{alt}{$\frac{1}{[\sigma^{(n_i)}]^{2}K_i}$}\psfrag{sigma}{$\frac{a_{i}-y(P_i)}{\sigma_{i}^{(n_{i})}
(\nu_i(0))}$}\psfrag{t}{$\leq k_i\lambda^{(n_i)}$}
\psfrag{madre}{secondary}
\psfrag{hija}{primary}\psfrag{bi}{$a_{i}$}\psfrag{nme}{$\nu_i^-$}
\begin{center}\includegraphics [scale=.25]{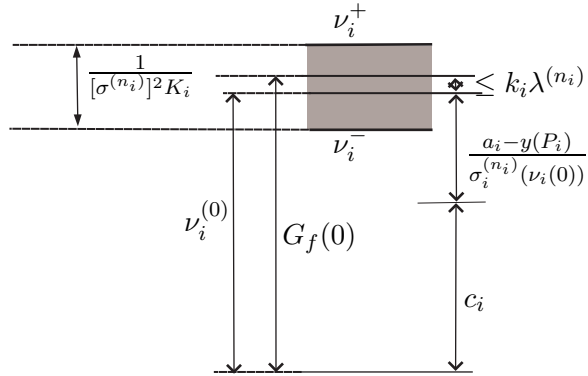}
\caption{\label{parametro} Determination of parameters for the sinks}\end{center}
\end{figure}

{\em Proof:} See figure \ref{parametro}. The equality
(\ref{lemabasicoigualdad}) in Lemma \ref{LemaBasico} which defines
the implicit function $\nu_i(\widehat \mu) = G_f(\widehat \mu )$ has
now the following expression, due to Remark
\ref{constantesenespaciofuncional}:
 \begin{equation}
 \begin{array}{rcl}
 \label{lemabasicoigualdad2}
\nu_i& =&c_i+ (a_i- y(P_i)) \cdot  [\sigma_{i}^{(n_i)}(\nu_i)]^{-1}
+\\ && + \widehat \mu  \cdot {\displaystyle
\frac{[\sigma_{i}^{(n_i)}(\nu_i)]^{-2}}{\beta_i(f)}}
  - \mbox{dist}^s(q_i,P_i) \cdot  \gamma_i (f)\cdot  [\lambda_{i}^{(n_i)}(b_i(f))]\\
 \end{array}
 \end{equation}

If $\widehat \mu = 0$ the equation (\ref{lemabasicoigualdad2})
depends on $f \in {\cal N}_1$ only because its last term does, and
defines $G_f(0) $.

On the other hand, the equation (\ref{lemabasico2igualdad}), which
is independent on $f \in {\cal N}_1$, defines
$\nu_i^{(0)}$.

Subtracting (\ref{lemabasico2igualdad}) and
(\ref{lemabasicoigualdad2}) with $\widehat{\mu}=0$ and applying the
Lagrange Theorem we obtain:
$$|G_f(0) - \nu_i^{(0)} | = \frac{|\mbox{dist}^s(q_i,P_i)  \gamma_i (f)
[\lambda_{i}^{(n_i)}(b_i(f))]|}{1+(a_i -
y(P_i))[\sigma_i^{(n_i)}]^{-2}[\sum_j
\sigma_{i,j}^{(n_i-1)}(Y)]\sigma'_{i,j}(Y)}$$ where $Y$ is an
intermediate value between $\nu_i(0)$ and $G_f(0)$. Using inequality
(\ref{valorespropios2}) and arguing as in Lemmas \ref{lemaBeta} and
\ref{nuevo2} we obtain:
$$|G_f(0) - \nu_i^{(0)} | \leq k_i [\max_{(x,y) \in \overline{V}_i} \lambda(x)]^{n_i}$$

%

We now compute its derivarive $G_f'$ respect to $\widehat \mu$ in
any point where $G_f$ is defined:
$$G_f ' = \frac{{{\displaystyle {\left[\sigma_{i}^{(n_i)}\right]^{-2}}}} }
{ {{\beta_i(f)} \cdot \left ( 1 +(a_i- y(P_i)) \cdot
\left[\sigma_{i}^{(n_i)}\right] ^{-2}\sum_j
\left[\sigma_{i,j}^{(n_i-1)}\right] \sigma_{i,j}'\right ) +2
\widehat \mu  \cdot {\displaystyle \left[\sigma_{i}^{(n_i)}\right]
^{-3}\sum_j
\left[\sigma_{i,j}^{(n_i-1)}\right] \sigma_{i,j}'{}}} }$$ 

Arguing as in Lemmas \ref{lemaBeta} and \ref{nuevo2} we conclude that it is uniformly bounded for $f \in {\cal N}_1$.

Finally, using  Lemma \ref{lemaNuevo}, the inequalities
(\ref{valorespropios2}) and recalling that  $n_i$ is large enough,
we deduce the bounds of $|G'_f(\widehat{\mu})|$ in the thesis.
$\;\;\Box$

\vspace{.2cm}

Combining the results in Proposition \ref{LemaBasico} and Lemma
\ref{Lemabasico2} we obtain the following:
\begin{lema}
\label{Lemabasico3} For each $i \geq 1$  and for each  sufficiently large $n_i$  there exist
constants  $ \nu_i^- < \nu_i^{0} < \nu_i^+$, independent of $f\in \mathcal{N}_1$,
such that:
$$|\nu_i^{+} - \nu_i^{-}| = \frac{1}{8(1+\varepsilon)  \beta_i(f_1)}
(\max_{y \in \overline{V}_i} \{\sigma  (y)\} )^{-2n_i}$$
$$ \nu_i^{(0) } = c_i + (a_i - y(P_i)) [ \sigma_{i}^{(n_i)}(\nu_i(0) ]^{-1}$$

and, if
\begin{equation}
\label{lemabasicoigualdad3} f \in {\cal N}_1: \; \; \;  \nu_i(f) \in
(\nu_i^-, \nu_i^+)
\end{equation} then $ f $ exhibits a sink $s_i (f)$
.  Even more, if $f, g \in {\cal N}_1$ are arc connected in ${\cal N}_1$ and all the connecting arc verifies the condition \em (\ref{lemabasicoigualdad3}), \em then the sink  $s_i(g)$ is the   continuation of the sink $s_i(f)$.
\end{lema}
{\em Proof: } See figure \ref{parametro}. Consider the $C^1$ real function $ \nu_i(f) = G_f(\widehat \mu)$ of real variable $\widehat \mu$ defined as in Lemma \ref{Lemabasico2}. It is strictly monotone because its first derivative is never zero. Applying the Lagrange Theorem to $G_f$, and the lower bound of its derivative given in Lemma \ref{Lemabasico2}, we deduce that the images by $G_f$ of the intervals  $ (-3/8,0] $ and $[0, 1/8)$ are  two intervals of length:
\begin{equation} \label{aux}
|G_f(1/8) - G_f(0)| = \frac{|G_f' (\widehat \mu _f^{(1)})|}{8}
> \frac{1}{8(1+\varepsilon)  \beta_i(f_1)}    (\max_{y \in \overline{V}_i} \{\sigma  (y)\} )^{-2n_i} = {C_i} >0
\end{equation}
$$|G_f(0) - G_f(-3/8)| =  \frac{3|G_f' (\widehat \mu _f^{(2)})|}{8}
> \frac{3}{8(1+\varepsilon)  \beta_i(f_1)}     (\max_{y \in \overline{V}_i} \{\sigma  (y)\} )^{-2n_i} = {3 C_i} >0 $$
where $C_i$ is a constant independent of $f \in {\cal N}_1$, but depending on $n_i$.

On the other hand, due to Lemma \ref{Lemabasico2}, there exists a
real number $\nu_i^{(0)}$, independent of $f \in {\cal N}_1$ such
that $$|\nu_i^{(0)} - G_f(0)| \leq K_i [\max_{(x,y)\in
\overline{V}_i}\lambda(x)] ^{n_i}, \;\; G_f(0) \in G_f((-3/8,
1/8))$$

As the horseshoe is strongly dissipative, $\lambda^{n_i} (x) \ll
\sigma ^{-2n_i}(y)$ if $n_i$ is large enough, for all points of
the rectangle $\overline{V}_i$, in particular, for the point $(x,y)$
where $\max \lambda(x)$ and $\max \sigma (y)$ are obtained. So we
can assume that $|\nu_i^{(0)} - G_f(0)| < C_i/2$ and then
$[\nu_i^{(0)} - C_i/2, \nu_i^{(0)} + C_i/2] \subset G_f((-3/8, 1/8))
$.

We define $\nu_i^{+} = \nu_i^{(0)} + C_i/2 $ and $  \nu_i^{-} =
\nu_i^{(0)} - C_i/2$.  Both values are independent of $f \in {\cal
N}_1$ and included in the image by $G_f$ of the interval $(-3/8,
1/8)$.

We use the definition of the
constant $C_i$ in the Equality (\ref{aux}) to get the exact value of $|\nu_i^+ - \nu_i^-|$.

The real function $G_f (\widehat{\mu})$ is strictly monotone.
Therefore, given $\nu_i(f) \in (\nu_i^{-}, \nu_i^+)$ there exists a
single value of $\widehat \mu \in (-3/8, 1/8)$ such that $\nu_i(f)=
G_f(\widehat \mu)$. Therefore, the hypothesis of   Proposition
\ref{LemaBasico} is fulfilled, and so  its thesis about the
existence and   continuation of the sink $s_i$ is verified.
$\;\;\Box$

\begin{prop}
\label{proposicion} The infinite dimensional  arc-connected manifold
${\cal M}_K \subset {\cal N}_1$ in Lemma \em \ref{variedadCod1}, \em
(which has codimension one in ${\cal N}_1$) verifies, for each $i
\geq 1$, the following properties:

(a) If $f\in {\cal M}_{K}$ with the constant $K= k(f_i \circ
f_1^{-1})$, then $f$ exhibits a homoclinic tangency at the point
$q_i \in V$ of the saddle $P_i \in \Lambda$.

(b) If $f\in {\cal M}_{K}$ with the constant $K \in k(f_i \circ
f_1^{-1})- c_i + (\nu_i^-, \nu_i^+) $, where $\nu_i^- < \nu_i^+$ are
defined as in Lemma \em \ref{Lemabasico3}, \em  then $f$ exhibits a
sink $s_i(f) \in V$. Even more, if $f, g \in {\cal M}_{K}$ then they
are isotopic and  the sink $s_i(g)$ is the   continuation of the
sink $s_i(f)$.
\end{prop}
{\em Proof: }

As proved in Lemma \ref{variedadCod1}, the manifolds ${\cal M}_K$
with $K$ constant are codimension one submanifolds of ${\cal N}_1$.

First choose, for each $i \geq 1$, a fixed $f_i \in {\cal N}_1$ such
that $f_i$ exhibits a homoclinic quadratic tangency at the point
$q_i \in V$. Due to the definition of $\mu_i(f)$ in subsection
\ref{paragrafo2}, such $f_i$ verifies $\mu_i(f_i) = 0$. Recall
Equalities (\ref{muconk}) and note that  for all $f \in {\cal N}_1$:
$ \mu_i(f) = 0$ if and only if $f$ exhibits a homoclinic
tangency at the point $q_i$.  This condition is fulfilled if and
only if $k(f \circ f_1^{-1}) = k(f_i \circ f_1^{-1})$, which proves
part a).

To prove part b) argue similarly, using Equalities (\ref{muconk})
with $\nu_i(f) \in (\nu_i^-, \nu_i^+)$, and applying Lemma
\ref{Lemabasico3}.  $\;\;\Box$

\vspace{.2cm}

\section{Proof of Theorem \ref{teoespesura}.} \label{seccionprueba1}

For the given diffeomorphism  $f_1$ as in the hypothesis of Theorem
\ref{teoespesura}, we shall work along the infinite dimensional
manifold ${\cal N}_1 \ni f_1$ of $C^3$ diffeomorphisms, defined in
Section \ref{seccionN}, by conditions (\ref{ambienteRestringido})
and (\ref{ambienteRestringido2}).

Let us construct, as in Newhouse-Robinson theorem (see
\cite{newhouse} and \cite{robinson}), a sequence of sinks $s_i$
which are produced along a monoparametric family of diffeomorphisms
(which we will call \em primary family\em), which generically
unfolds a sequence of homoclinic quadratic tangencies $q_i$. By
induction, each sink $s_{i+1}$ shall be produced while the $i$ sinks
that were previously generated, still survive. The key to get this
result is the persistence of tangencies of the Theorem of Newhouse
in \cite{newhouse}. For a seek of completeness we reproduce here the
details of the inductive proof of Newhouse-Robinson Theorem to
obtain infinitely many simultaneous sinks. We improve the argument,
adding the conclusions of our previous sections, to obtain also the
simultaneous continuation of the infinitely many sinks.

\begin{defi} \em \label{primarysecondary}

Given a small enough real number $\delta >0$ we fix a one -parameter
family $\{\widetilde{f}_{t}\}_{t \in (- \varepsilon, \varepsilon)}
\in {\cal N}_1 $, called the \em primary family, \em such that:

$\widetilde{f}_t= \xi_t \circ f_1$ where $\xi_t$ verifies the
conditions (\ref{ambienteRestringido2}) and besides:

$\| \xi_t - \mbox{id}\|_{C^3} < \delta \forall \, t \in
(-\varepsilon,\varepsilon)$

$\xi_0= Id $, and thus $\widetilde{f}_0 =f_1$ and $ \; k(\xi_0) = 0$

$k(\xi_t) = t$

We call \lq \lq secondary family" $g_t$ to any other one-parameter
family $\{g_{t}\}_{t \in (- \varepsilon, \varepsilon)} \in {\cal
N}_1$ such that  $k(g_t \circ f_1^{-1}) = k(\widetilde{f}_t \circ
f_1^{-1}) = t$ for all $t \in (- \varepsilon,  \varepsilon)$.
\end{defi}

We observe that the primary family is transversal to the manifolds
${\cal M}_K$ defined in Lemma \ref{variedadCod1}, for all $K \in
(-\varepsilon, +\varepsilon)$, and that it unfolds generically any
homoclinic tangency $q_i$ produced in the line of tangencies
$L_{\widetilde{f}_t}$, in particular the given homoclinic tangency
$q_1 \in L_{f_1}$. After the density of tangencies (see \cite{newhouse}) the
hypothesis of large thickness of the stable and unstable Cantor sets
$K^s$ and $K^u$ of $\Lambda $ along the line of tangencies $L_f$,
assumed in Theorem \ref{teoespesura}, implies the following:
%

 \begin{obs}
\label{persistencetangencies} If $\varepsilon >0$ is small enough
then there exists a dense set of parameter values $t \in
(-\varepsilon,\varepsilon)$ in the primary family, such that the
diffeomorphism $\widetilde{f}_t$ exhibits a homoclinic tangency in
some point $q $ in the line of tangencies $L_{\widetilde{f}_t}$ of
some saddle periodic point $P \in \Lambda$.
 \end{obs}

%

We take the system of coordinates in the open neighborhood $U
\supset \Lambda$, as defined in Section \ref{seccionCoordenadas} and
such that the saddle $P_1 = (0,0)$. We shall choose the tangency
$q_1 \in U$ such that $q_1= (b_1, 0)$ (i.e. in the connected local
stable leaf of $P_1$), and then the small neighborhood $q_1 \in V
\subset U$ (see Remark \ref{defincionV}). We use the notation of
Section \ref{seccionN}:

We have $t= 0$, $ k(\xi_0) = k(id) =  0  $ and $\mu _1(f_1) = 0 $, $
\nu_1 (f_1) = c_1 + \mu_1 (f_1) = 0$ and there is a tangency at the
point $q_1 = (b_1, 0) \in L_1= L(f_1)$.

Applying Lemma \ref{Lemabasico3} we choose $n_1$ sufficiently large
so the fixed numbers $\nu_1^-$ and $\nu_1^+$ verify $|\nu_1^{\pm}|<
\varepsilon$. Applying Proposition \ref{proposicion}, if  $t  =
k(\widetilde{f}_t \circ f_1^{-1})  \in (\nu_1^-, \nu_1^+) $ then
$\widetilde{f}_t$ exhibits a sink $s_1(\widetilde{f}_t)$ and for all
secondary diffeomorphism $g \in {\cal N}_1$ such that $k(g \circ
f_1^{-1}) = t$ there exists the continuation $s_1 (g)$ of the sink
$s_1(\widetilde{f}_t)$.

We now argue by induction in the number of simultaneous sinks
exhibited by $\widetilde{f}_t$ in the primary family:

Let us suppose that there exist parameter values $-\varepsilon <
t_i^-< t_i^+< \varepsilon  $ such that for all $ t \in
(t_i^-,t_i^+)$ the diffeomorphism $\widetilde{f}_t$ of the primary
family exhibits $i$ sinks $s_1(\widetilde{f}_t),
s_2(\widetilde{f}_t), \ldots, s_i(\widetilde{f}_t) \in V$, and  any
secondary diffeomorphism $g_t$  exhibits  the   continuations
$s_1(g_t), s_2(g_t), \ldots, s_i(g_t) \in V $ of those $i$ sinks.

We shall construct an  interval with non void interior $[t_{i+1}^- ,
t_{i+1}]^+  \in (t_i^-,t_i^+)$ such that if $t \in (t_{i+1}^- ,
t_{i+1}^+) $ then the diffeomorphisms $\widetilde{f}_t$ of the
primary family exhibit a new sink $s_{i+1}(\widetilde{f}_t) \in V$,
and besides all secondary diffeomorphisms $g_t$ with $t \in
(t_{i+1}^- , t_{i+1}^+)$ exhibit the   continuation $s_{i+1}(g_t)$
of $s_{i+1}(\widetilde{f}_t)$.

After the density of tangencies (see Remark \ref{persistencetangencies}), there exists  a parameter value \begin{equation}
\label{puntomedio} t_{i+1} \in \left ( t_i^- + \frac{ t_{i}^+-t_{i}^-}{4},t_i^- + \frac{3(t_{i}^+-t_{i}^-)}{4} \right )
\end{equation}
such that some
periodic saddle point $P_{i+1}$ in the horseshoe $\Lambda$, has a
homoclinic tangency at $q_{i+1} = (b_{i+1}, c_{i+1}) \in L (\widetilde{f}_{t_{i+1}})$.

Unfolding the tangency of $f_{{i+1}}=\widetilde{f}_{t_{i+1}}$ when
moving along the primary family, we will create a new sink $s_{i+1}$
in such a way that $t$ is still in the interval $(t_{i}^-, t_i^+)$
where the $i$ previous sinks still persist.

To construct the new sink $s_{i+1}$ and a parameter interval inside
$(t_{i}^-, t_i^+)$  in which this new sink is exhibited, we argue as
follows:

If $t= t_{i+1}$ verifies condition (\ref{puntomedio}),   there
exists a homoclinic tangency at $q_{i+1} = (b_{i+1}, c_{i+1})$ of a
saddle $P_{i+1} \in \Lambda$. Therefore the height of the parabolic
unstable arc of $P_{i+1}$ is $\mu_{i+1}(f_{{i+1}})= 0 , \; \nu_{i+1}
= c_{i+1}$ (see Figure \ref{gr10}). On the other hand, by Definition
\ref{primarysecondary} of the primary family, we have: $ t_{i+1} =
k(f_{{i+1}} \circ f_1^{-1}) $.

Applying Lemma \ref{Lemabasico3}, we shall find the parameter values
$t_{i+1}^-< t_{i+1}^+ $ and the diffeomorphisms $\widetilde{f}_{t}$
in the primary family such that
$\nu_{i+1}(\widetilde{f}_t) \in [\nu_{i+1}^-, \nu_{i+1}^+]$
if $t \in [t_{i+1}^-, t_{i+1}^+] \subset (t_{i}^-, t_{i}^+)$.
Also from Lemma \ref{Lemabasico3}, we can choose a sufficiently
large $n_{i+1}$ so that $|\nu_{i+1}^{\pm} - c_{i+1}| < (t_i^+- t_i^-)/8$.

After equalities (\ref{muconk}), and recalling that $t =
k(\widetilde{f}_t \circ f_1^{-1})$ we obtain:
$\nu_{i+1}(\widetilde{f}_t) - c_{i+1} = t - t_{i+1}$. So, if $t \in
t_{i+1} - c_{i+1} + (\nu_{i+1}^-, \nu_{i+1}^+)$, then
$\nu_{i+1}(\widetilde{f}_t) \in [\nu_{i+1}^-, \nu_{i+1}^+] $, and
applying Lemma \ref{Lemabasico3}, the map $\widetilde{f}_t$ will
exhibit a new sink $s_{i+1}(\widetilde{f}_t)$.

Applying Proposition \ref{proposicion}, for any secondary family
$\{g_t\}_t$, if $t = k(g \circ f_1^{-1}) \in t_{i+1}- c_{i+1} +
[\nu_{i+1}^-, \nu_{i+1}^+]$  then the map $g_t$ will exhibit the
continuation $s_{i+1}(g_t)$ of the sink $s_{i+1}(\widetilde{f}_t)$.

We then define $$t_{i+1}^{\pm} = t_{i+1}- c_{i+1} + \nu_{i+1}^{\pm}\in
(t_{i+1} -(t_i^+- t_i^-)/8, t_{i+1} +(t_i^+- t_i^-)/8) $$ From
condition (\ref{puntomedio}), the equality above implies
$[t_{i+1}^-, t_{i+1}^+] \subset (t_i^-, t_i^+)$ as wanted. We
conclude that for $t \in (t_{i+1}^-, t_{i+1}^+)$ the diffeomorphisms
$\widetilde{f}_t$  of the primary family  exhibit the sink
$s_{i+1}(\widetilde{f}_t)$ and the secondary diffeomorphisms $g_t$
exhibit the   continuation $s_{i+1}(g_t)$ of that sink.

We observe that we can make this new sink $s_{i+1}$ of arbitrary
period $n_{i+1}+ N_{i+1}$, provided that it shall be  large enough,
because in Lemma \ref{Lemabasico3} we can arbitrarily choose the
natural number $n_{i+1}$, from a minimum value.  If we choose $n_{i+1}$ such
that $n_{i+1} + N_{i+1}$ is not a multiple of the periods of the previous
$i$ sinks $s_1, s_2, \ldots, s_i$, then  the sink $s_{i+1}$ shall be
necessarily a new one.

Finally, taking $t_{\infty} = \bigcap _{i=1}^{\infty} [t_{i}^-,
t_{i}^+]$, as we constructed each compact interval in the interior
of the previous one, the real value $t_{\infty}$ is in the interior
of all intervals, and thus, by construction, there exists the sink
$\{s_{i}(g)\}$ for all $i \geq 1$ and for all diffeomorphism $g \in
{\cal M}_{t_{\infty}}$, being each sink $s_i(g)$  the   continuation
of the respective sink $s_i(\widetilde{f}_{t_{\infty}})$ exhibited
by the diffeomorphism $\widetilde{f}_{t_{\infty}}$ in the primary
family. This  ends the proof of Theorem \ref{teoespesura}.
$\;\;\Box$

\section{Conclusion of the main results.} \label{seccionpruebateo}

{\em End of the Proof of Theorem } \ref{teo}:

Due to Newhouse-Robinson Theorem \ref{newhouse}, we find a
diffeomorphism $f_1 \in {\cal N}$, arbitrarily near the given $f_0$,
such that $f_1$ has a horseshoe $\Lambda$, which is strongly
dissipative and fulfils the condition of large thickness, and
besides there is an homoclinic quadratic tangency $q_1$ of a saddle
$P_1 \in \Lambda$. These last assertions are the hypothesis of
Theorem \ref{teoespesura} which we have already proved in the last
section, ending the proof of Theorem \ref{teo}. $\;\; \Box$

\vspace{1cm}

{\em Proof of part A) of Theorem } \ref{teoprevio} (see Figure
\ref{tg}): The given one-parameter family  $\{\widetilde{f}_t\}_{t
\in (-\varepsilon,\varepsilon)}$ generically unfolds the quadratic
tangency at $q_0$ of a saddle point $P_0$. The generic unfolding is
defined by the condition. $v= d \mu _0 (\widetilde{f}_t)/dt \neq 0$
for all $t \in (-\varepsilon, \varepsilon)  $, where $\mu _0
(\widetilde{f}_t)$ is the height of the parabolic arc in the
unstable leaf of the saddle $P_0$ respect to the local stable arc of
$P_0$, to which it is tangent when $t= 0$ (i.e.
$\mu_0(\widetilde{f}_t)|_{t=0}=0$).

After the Theorem of Newhouse-Robinson (revisited in Theorem
\ref{newhouse}), there exists an interval $I \subset
(-\varepsilon,\varepsilon)$ such that for all $t \in I$ there is a
horseshoe $\Lambda$ fulfilling the large thickness condition (see
Definition \ref{espesurasgrandes}) which is strongly dissipative
(see Lemma \ref{disipacionfuerte}). Even more, Newhouse-Robinson
Theorem asserts that there exist a dense set $H \subset I$ of
parameter values, such that $\widetilde{f}_t$ exhibits some
homoclinic tangency,  for all $t \in H$. Let us choose some $t_1 \in
H$, such that $\widetilde{f}_{t_1}=f_1$ exhibits such homoclinic
tangency at a point $q_1 \in V$ (see remark \ref{defincionV}) of a
saddle $P_1 \in \Lambda$.

Consider for this $f_1$ a small neighborhood ${\cal N}$ such that,
for all $f \in {\cal N}$ (in particular for all $\widetilde{f}_t$ with $t $ near
$t_1$), there exist the real numbers $a_1(f), b_1(f), c_1(f), \alpha_1(f),
\beta_1(f), \gamma_1(f), \nu_1(f), \mu_1(f)$ defined in Section
\ref{seccionCoordenadas},  Figure \ref{tg} and Equations
(\ref{alphabetagamma}).


Observe that $|d \mu _1(\widetilde{f}_t) /dt|\geq | v |
/2 \; \neq 0$. Suppose $d \mu _1(\widetilde{f}_t) /dt >0$, then $\mu_1(\widetilde{f}_t)$ is a strictly
increasing diffeomorphic function which is zero for $t=t_1$, and
there exists $\delta_1 >0$ such that $\mu_1 (\widetilde{f}_t) \in
(-\delta_1, \delta_1) \; \Rightarrow \; \; t \in I$.

We will repeat the well known argument of Newhouse, improving it to
get the eigenvalues of the sinks as small as wanted:

First, we shall construct some interval $(t_1^-, t_1^+) \subset I$
of the parameter values $t$ for which $\widetilde{f}_t$ exhibits a
sink in $V$, whose eigenvalues have modulus smaller than the given
number $0<\rho <1$.

\vspace{.2cm}

Consider, for each $\widetilde{f}_t$, the first equation
(\ref{reparam3})  giving $\mu_1(\widetilde{f}_t)$ diffeomorphically
as a function of the new parameter $\widehat \mu$, for each fixed
$n=n_1 \geq 1$. If $n_1$ is large enough, there exist  $ -\delta_1
<\mu_1^-< \mu_1^+< \delta_1  $ such that if $\mu_1 (\widetilde{f}_t)
\in (\mu_1^-, \mu_1^+)$ then $ \widehat \mu \in (k(\rho)^-,
k(\rho)^+)$, defined in Proposition \ref{LemaBasico}. Therefore, the
thesis of this proposition implies that $\widetilde{f}_t$ has a sink in
$V$ whose eigenvalues have modulus smaller than $\rho$. Considering
that the real function $\mu_1(\widetilde{f}_t)$ depends
diffeomorphically on $t$, the preimage by $\mu_1(f_{t})$ of the
interval $(\mu_1^-, \mu_1^+)$ is an interval $(t_1^-, t_1^+) \subset I$.
By construction, if $t \in (t_1^-, t_1^+)$ then $\widetilde{f}_t $
exhibits a sink $s_1$ in $V$ whose eigenvalues have modulus smaller
than $\rho$.

Now, by induction, suppose that there is an open interval $(t_i^-, t_i^+)
\subset I$ such that, if $t \in (t_i^-, t_i^+) $, then
$\widetilde{f}_t $ exhibits $i$ simultaneous different orbits of the
sinks $s_1, s_2, \ldots, s_i$ in $V$, whose eigenvalues have all
modulus smaller than $\rho$. As the set $H \subset I$, where the
homoclinic tangencies are produced, is dense in $I$, we can choose
$t_{i+1} \in (t_i^- + (1/4) (t_i^+ - t_i^-), \; \; t_i ^- + (3/4)
(t_i^+ - t_i^-))$ such that $\widetilde{f}_{t_{i}}$ exhibits a
homoclinic point $q_{i+1} \in V$ of a saddle $P_{i+1}$. As above,
the function $\mu_{i+1}(f_t)$ is an increasing  diffeomorphism from
the interval $(t^-_{i}, t^+_i)$, to an interval of the real variable
$\mu_{i+1}(\widetilde{f}_t)$, such that
$\mu_{i+1}(\widetilde{f}_{t_{i+1}}) = 0$. Therefore, there exists
$\delta_{n+1} >0$ such that, if $\mu_{i+1} (\widetilde{f}_t) \in
(-\delta_{n+1}, \delta _{n+1})$, then $|t - t_{i+1}| < 1/8(t_i^+-
t_i^-) $.

Arguing as in the first step, if $n_{i+1}$ is large enough there
exist $ -\delta_{n+1} <\mu_{n+1}^-< \mu_{n+1}^+< \delta_{n+1} $ such
that if $\mu_{n+1} (\widetilde{f}_t) \in (\mu_{n+1}^-, \mu_{n+1}^+)$
then $ \widehat \mu \in (k(\rho)^-, k(\rho)^+)$, defined in Proposition
\ref{LemaBasico}. Therefore, the thesis of this proposition implies that
$\widetilde{f}_t$ has a sink $s_{i+1}$ in $V$ whose eigenvalues have
modulus smaller than $\rho$. This new sink is different from the $i$
sinks that were previously constructed, provided one can choose any
integer number $n_{i+1}$ large enough, so one can get the period of
the new sink larger and not a multiple, of the periods of the $i$
sinks that were previously constructed.

Considering that the real function $\mu_{i+1}(\widetilde{f}_t)$
depends diffeomorphically on $t$, the preimage by $\mu_{i+1}(f_{t})$
of the interval $(\mu_{i+1}^-, \mu_{i+1}^+)$ is an interval
$$[t_{i+1}^-, t_{i+1}^+] \in t_{i+1} + [-(t_{i}^+ -t_i^-)/8, (t_{i}^+.
- t_i^-)/8] \subset (t_i^-, t_i^+)$$ By construction, if $t \in
(t_{i+1}^-, t_{i+1}^+)$ then $\widetilde{f}_t $ exhibits
simultaneously $i+1$  sinks in $V$, whose eigenvalues have modulus
smaller than $\rho$.

Finally define $g_0 = \widetilde{f}_{t_{\infty}}$ where $t_{\infty}
\in \bigcap _{i=1}^{\infty} [t_{i}^-, t_i^+]$.  By construction, the
set $J$ of such values where the inifinitely many sinks exist is
dense in $I$. $\;\;\Box$

\vspace{.3cm}

\begin{obs} \em \label{remarkmsubi}
We observe that in in the proof of part (A) of Theorem
\ref{teoprevio}, the construction of the map $g_0 = f_{t_{\infty}}$,
which exhibits infinitely many sinks in $V$, allows us to choose,
for each $i \geq 1$, any integer $n_i$ provided it is large enough.
So, we can obtain the same thesis if, besides, we ask $n_i \geq
m_i$, where $m_i \rightarrow +\infty$ is any previously specified
sequence of integer numbers.
\end{obs}

\vspace{.3cm}

{\em Proof of part B) of Theorem } \ref{teoprevio}:
We will show that
there exists a sequence $m_i \rightarrow + \infty$ such that, if
$g_0$ is constructed as in the proof of part (A) and besides
verifying $n_i > m_i$ for all $i \geq 1$, then the thesis (B) of
Theorem \ref{teoprevio} holds for this $g_0$.

Choose $\delta >0$ small enough (to be fixed at the end of the
proof) and define the following manifold ${\cal N}_1 \subset
\mbox{Diff}^3(M)$, which is $\delta- C^3$ near $g_0$:
$${\cal N}_1= \{g \in \mbox{Diff}^3(M): g= \xi \circ g_0    \}$$
where $\xi \in \mbox{Diff}^3(M)$ is such that
$\|\xi-\mbox{id}\|_{C^3}<\delta$ and besides it verifies conditions
(\ref{ambienteRestringido2}), replacing $g_0$ instead of $f_1$, in a
small fixed neighborhood $V$ of the line of tangencies $L_0 = L(g_0)
$ (instead of the line of tangencies $L_1 = L(f_1)$).

As in equalities (\ref{muconk}), we now  have $$k(g \circ g_0^{-1}) = \mu_i(g) - \mu_i(g_0) = \nu_i(g) - \nu_i(g_0) \; \; \forall \; i \geq 1$$

Consider the set ${\cal {M}}$,  as follows:
$${\cal M} = \{g \in {\cal N}_1: \; k(g \circ g_0^{-1}) = 0\}$$
As seen in Lemma \ref{variedadCod1}, ${\cal M}$ is an infinite
dimensional, arc connected manifold, with codimension one in ${\cal
{N}}_1$. By construction:
\begin{equation}
\label{muconk2}g \in {\cal M}  \; \; \Rightarrow \; \; \mu_i(g) = \mu_i(g_0) , \; \;  \nu_i(g) = \nu_i(g_0) \; \; \forall \; i \geq 1
\end{equation}

Let us prove that, if $g \in {\cal M}$, then the infinitely many
sinks $s_i(g_0)$  have continuation sinks $s_i(g)$.

We apply, for $g$ and
$g_0$, the respective changes of variables and parameter given by
Equations (\ref{reparam}) and (\ref{reparam3}). We
recall from the proof of part (A) of Theorem \ref{teoprevio}, that
$\mu_i (g_0)$ was constructed such that: $$\widehat{\mu} (g_0) \in
(k^-(\rho), k^+(\rho))$$ where $$ -\frac{3}{8} < k^{-}(\rho) < 0 < k
^{+}(\rho) < \frac{1}{8}$$ are the numbers defined in Proposition \ref{LemaBasico}.

By contradiction, if $g \in {\cal M} \subset {\cal N}_1$ did not
have the continuation of the sink $s_i(g_0)$, then,  due to
Proposition \ref{D} we obtain
$\widehat \mu (g) \not \in (-3/8, 1/8)$.
In fact, we note that  $g_0$ and $g$ are not isotopic by any
one-parameter family of diffeomorhphisms $\{g_t\}_{0 \leq t \leq
1}$ such that
$\widehat \mu (g_t) \in (-3/8, 1/8) \; \forall t \in [0,1]$. So
considering in particular some one-parameter family in ${\cal M}$,
there would exist 
an interval $[t_0,t_1]$ such that $\widehat{\mu}(g_{t_0}) \in \{
k^-(\rho),k^+(\rho)\}$, $\widehat{\mu}(g_{t_1}) \in \{ -3/8,1/8\}$
and $\forall \, t \in [t_0,t_1]$, $\widehat{\mu}(g(t)) \neq 0$.

Therefore:
$$\left|\frac{\widehat \mu (g_{t_0})-\widehat \mu (g_{t_1}) } {\widehat \mu (g_{t_1})}\right|
\geq \frac{\min \{|1/8 - k^+(\rho)|, | k^-(\rho)-3/8|\}}{3/8} = \eta (\rho) = \eta >0$$
\begin{equation} \label{ponernumero}
0< \frac{\widehat \mu (g_{t_0})}{\widehat \mu (g_{t_1})}
 \not \in (1 - \eta, 1 + \eta)\end{equation}


From equation (\ref{muconk2}), we observe that $\mu_i(g_{t_1}) =
\mu_i(g_{t_0})$. Then, taking into account (\ref{reparam3}) in the
points $ y = \nu_i(g_{t_1}) = \nu_i(g_{t_0})$, $ x= b_i(g_{t_1})
\mbox{ or } x= b_i(g_{t_0})$ and subtracting:
\begin{equation}
\label{diferenciamutechos} 0 = \frac{\widehat \mu(g_{t_1})}{\beta_i(g_{t_1})}
(\sigma_{i}^{(n_i)})^{-2} - \frac{\widehat \mu (g_{t_0})}{\beta_i(g_{t_0})}
(\sigma_{i}^{(n_i)})^{-2} + R_i(g_{t_1}, g_{t_0}) \end{equation} where $R_i(g_{t_1},
g_{t_0})$ is obtained as the difference of the two terms (one computed
for $g_{t_1}$ and the other for $g_{t_0}$) in Equations (\ref{reparam3}), that
have the factor $\lambda_{i}^{(n_i)}(x)$. By the strong dissipative
condition we have $\lambda_{i}^{(n_i)} \ll (\sigma_{i}^{(n_i)})^{-2}$ if $n_i
$ is large enough. Besides, the other terms or coefficients
depending continuously on $g$ in the Equation (\ref{reparam3}), are
upper and lower bounded from zero, due to Lemma \ref{lemaBeta}.
Therefore given $0 <\varepsilon $ there exists $m_i$ large enough such that if $n_i\geq m_i$ then 
$$|R_i(g_{t_1}, g_{t_0})|\leq k_i \lambda_i^{(n_i)} \leq \frac{\varepsilon }{8 \sup_{g \in {\cal N}} |\beta_i(g)|} \; \; (\sigma_{i}^{n_i})^{-2} \leq \frac{\varepsilon |\widehat \mu(g_{t_1})|}{  |\beta_i(g_{t_1})|}\; \; (\sigma_{i}^{n_i})^{-2} $$
Substituting in \ref{diferenciamutechos} we obtain:
\begin{equation}\label{823} \frac{1-\varepsilon}{1+\varepsilon} \leq \frac{1}{1+\varepsilon}\left| \frac{\beta_i(g_{t_1})}{ \beta_i(g_{t_0})}\right|\leq \left| \frac{\widehat{\mu}(g_{t_1})}{\widehat{\mu}(g_{t_0})}\right| \leq \frac{1}{1-\varepsilon}\left| \frac{\beta_i(g_{t_1})}{ \beta_i(g_{t_0})}\right|\leq \frac{1+\varepsilon}{1-\varepsilon}\end{equation}
In the last inequalities we have used Lemma \ref{lemaNuevo}. Take
$\varepsilon >0$ such that
$$\left(\frac{1-\varepsilon}{1+\varepsilon},
\frac{1+\varepsilon}{1-\varepsilon}\right) \subset (1-\eta,1+\eta)$$
and then fix $\delta = \delta(\varepsilon)$ as in Lemma
\ref{lemaNuevo} to define $\mathcal{N}_1$. Therefore (\ref{823})
implies $$\left|
\frac{\widehat{\mu}(g_{t_1})}{\widehat{\mu}(g_{t_0})}\right|\in
(1-\eta,1+\eta)$$ contradicting (\ref{ponernumero}).   $\;\;\Box$

\vspace{.5cm}

{\em Proof of part C) of Theorem } \ref{teoprevio}:

Consider $g_0$ constructed as in the proof of the part (A) of
Theorem \ref{teoprevio}. The given one-parameter family $\{g_t\}_{t
\in (-\varepsilon,  \varepsilon)}$ is not  necessarily in the space
${\cal N}_1$ defined in that proof, but nevertheless it unfolds
generically the tangencies along $L(g_t)$, i.e.:
$$\left|\frac{d\mu_i(t)}{dt} \right| \geq v >0 \; \;
\forall i \geq 1, \; \; \forall t \in (-\varepsilon, \varepsilon)$$
where $\mu_i$ is defined in subsection \ref{paragrafo2}. We have:
$$|\mu_i(g_t) - \mu_i(g_0)| = |(\nu_i(g_t) - c_i(g_t)) -(\nu_i(g_0) - c_i(g_0))| \geq v |t| \; \; \forall \; i \geq 1$$

%
Consider the implicit function $\nu_i(f) \in
G^{(i)}_f(\widehat{\mu})$  verifying equation
(\ref{lemabasicoigualdad}) for any $f \in \mathcal{N}$. We can not
apply directly the thesis of Lemma \ref{Lemabasico2} because it is
valid only if $f \in \mathcal{N}_1$ and our diffeomorphism $g_t$
does not necessarily belong to the manifold $\mathcal{N}_1$.
Nevertheless we use equation (\ref{lemabasicoigualdad}), Lemma
\ref{lemaBeta} and similar arguments of those in the proof of Lemma
\ref{Lemabasico2} to obtain the following bounds for all $f \in
\mathcal{N}$ and for all $n_i$ large enough:
\begin{equation}
\label{824} |G_f^{(i)}(0) - c_i(f)| = \left|\frac{a_i(f) - y(P_i(f))}{\sigma_i^{(n_i)}(G_f^{(i)}(0))} - \mbox{dist}^s(q_i(f).P_i(f))\gamma_i(f)\lambda_i(f)^{(n_i)}(b_i(f))\right|<\frac1i
\end{equation}

\begin{equation}
\label{825} \left|\frac{dG_f^{(i)}(\widehat{\mu})}{d\widehat{\mu}}\right| \leq \frac{2}{K_i (\sigma_i^{(n_i)}(G_f^{(i)}(\mu)))^2}
\end{equation}

In particular we apply (\ref{824}) to $f=g_t$ and $f=g_0$ to obtain
\begin{equation}  \label{830} |(G ^{(i)}_{g_t}(0) - c_i(g_t)) - (G ^{(i)}_{g_0}(0)- c_i(g_0))|_{i \rightarrow +\infty } \rightarrow 0
 \end{equation}

Suppose that  $g_t$ exhibits infinitely many sinks $s_i(g_t)$ that are continuations from those of $g_0$.
Applying part (C) of Proposition \ref{LemaBasico} we deduce that there exists $\widehat{\mu} \in (-1,1)$ such that
\begin{equation}
\label{826} \nu_i(g_t) = G_{g_t}(\widehat{\mu})
\end{equation}
Combining (\ref{826}) with (\ref{825}) applied to  $f=g_t$  and
using the Lagrange Theorem:
$$\nu_i(g_t) \in G^{(i)}_{g_t} (0) + \frac{ [\sigma_{i}^{(n_i)}(g_t)(Y_i(t))]^{-2}}{K_i}
\left (-{2},  {2} \right )$$ But $\sigma_{i}^{(n_i)}(g_t)(Y_i(t))>
\widetilde{\sigma}^{n_i}$ with  $\widetilde{\sigma}>1$ and $n_i \to
\infty$ as fast as needed.

Recalling (\ref{830}):
$$ 0 < v , \; \; \; 0 \leq v |t|  \leq |(\nu_i(g_t) - c_i(g_t)) -
(\nu_i(g_0) - c_i(g_0))|_{i \rightarrow + \infty} \rightarrow 0 $$
and then $ |t| = 0$. $\;\; \Box$


\begin{thebibliography}{XXXXXX 2222}




\bibitem[HPS 1977]{hirshpughshub} \B{M. Hirsch, C. Pugh and M Shub}
{Invariant manifolds}
{Lecture Notes in Mathematics 583}{1977}

\bibitem[M 1973]{demelo} \A{de Melo}{Structural stability of diffeomorphisms on two-manifolds}
{Inventiones Math}{21}{1973}{233-246}

\bibitem[N 1974]{newhouse} \A{S. Newhouse}{Diffeomorphisms with infinitely many sinks}
{Topology}{13}{1974}{9-18}

\bibitem[P 2000]{palis} \A{J. Palis}{A global view of dynamics and a conjecture on the denseness of
finitude of attractors. \em G{\'e}om{\'e}trie complexe et syst\`{e}mes
dynamiques (Orsay, 1995)}{Ast{\'e}risque}{261}{2000}{335-347}


\bibitem[PT 1993]{palistakens} \B{J. Palis and F. Takens}
{Hyperbolicity and sensitive chaotic dynanics of homoclinic bifurcations}
{University Press, Cambridge}{1993}

\bibitem[R 1983]{robinson} \A{C. Robinson}{Bifurcation to infinitely many sinks}
{Comm Math Phys.}{90}{1983}{433-459}

\bibitem[YA 1983]{yorke} \A{J. A. Yorke and K. T. Alligood}
{Cascades of period doubling bifurcations: a prerequisite for horseshoes}{Bull AMS}{9}{1983}{319-322}




\end{thebibliography}
\end{document}